\documentclass[9pt]{article}
\usepackage[latin1]{inputenc}
\usepackage{amscd}
\usepackage{amsfonts}
\usepackage{amsgen}
\usepackage{amsmath}
\usepackage{amssymb}
\usepackage{amstext}
\usepackage{amsthm}
\usepackage{graphicx}
\usepackage{indentfirst}
\usepackage{latexsym}
\usepackage{makeidx}
\usepackage{mathrsfs}
\usepackage{epsfig}
\usepackage{wrapfig}
\usepackage[all,knot,arc,import,poly]{xy}
\usepackage[usenames]{color}
\usepackage{indentfirst}
\usepackage{mathrsfs}
\usepackage{epsfig}
\usepackage{wrapfig}
\usepackage[all,knot,arc,import,poly]{xy}
\input{xypic}

\newtheorem{Df}{Definition}[section]
\newtheorem{Teo}[Df]{Theorem}
\newtheorem{Prop}[Df]{Proposition}
\newtheorem{Lem}[Df]{Lemma}

\newtheorem{Obs}[Df]{Remark}
\newtheorem{Fat}[Df]{Fact}

\newtheorem{Cor}[Df]{Corollary}

\newtheorem{Ct}[Df]{}

\newcommand{\n}{\noindent}
\newcommand{\Dem}{\n{\em Proof:\;\;}}

\newcommand{\bc}{\begin{center}}
\newcommand{\ec}{\end{center}}

\newcommand{\N}{\mathbb{N}}

\newcommand{\vsete}{\vspace*{0.7cm}}
\newcommand{\vcinco}{\vspace*{0.5cm}}
\newcommand{\vtres}{\vspace*{0.3cm}}

\newcommand{\hem}{\hspace*{1em}}

\newlength{\dede}     

\newcommand{\cI}{\mbox{$\cal I$}}

\newcommand{\cK}{\mbox{$\cal K$}}
\newcommand{\cA}{\mbox{$\cal A$}}
\newcommand{\cQ}{\mbox{$\mathcal{Q}$}}
\newcommand{\cS}{\mbox{$\cal S$}}
\newcommand{\cT}{\mbox{$\mathcal{T}$}}
\def \rest {\restriction}
\def \ra {\rightarrow}

\def \lra {\longrightarrow }

\def \rat {\rightarrowtail }
\def \thra {\twoheadrightarrow }

\def \hookr {\hookrightarrow }

\newcommand{\Af}{\cA_f}
\newcommand{\Ss}{\cS_s}
\newcommand{\Sf}{\cS_f}
\newcommand{\Ls}{\mathcal{L}_s}
\newcommand{\Lf}{\mathcal{L}_f}
\newcommand{\Qf}{\cQ_f}
\newcommand{\Qc}{\cQ_f^c}
\def \beq { \begin{equation} }
\def \eeq { \end{equation} }

\def \rest {\restriction}

\def \sub {\subseteq}

\def \bu {\bullet}

\begin{document}
\title{Algebraizable Logics and a functorial encoding of its morphisms}
\author{H.L. Mariano, D.C. Pinto}

\maketitle

\begin{abstract}
The present work presents some results about the categorial relation between logics and its categories of structures. A (propositional, finitary) logic is a pair given by a signature and Tarskian consequence relation on its formula algebra. The logics are the objects in our categories of logics; the morphisms are certain signature morphisms that are translations between logics (\cite{AFLM1},\cite{AFLM2},\cite{AFLM3} \cite{FC}). Morphisms between algebraizable logics (\cite{BP}) are translations that preserves algebraizing pairs (\cite{MaMe}): they can be completely encoded by certain functors defined on the quasi-variety canonically associated to the algebraizable logics. This kind of results will be useful in the development of a categorial approach to the representation theory of general logics (\cite{MaPi1}, \cite{MaPi2}, \cite{AJMP}).
\end{abstract}

\section{Introduction}

The main motivation for the definition and development of categories of logics  is the combining logics methods. In the 1990's rise many methods of combinations of logics (\cite{CC3}). They appear in dual aspects: as processes of decomposition or analysis of logics (e.g., the "Possible Translation Semantics" of W. Carnielli, , \cite{Car}) or as processes  of  composition or synthesis of logics (e.g., the "Fibrings" of D. Gabbay, \cite{Ga}). The combining of logics is still a young topic in contemporary logic. Besides the pure philosophical interest of define mixed logic systems in which distinct operators obey logical relations  of different nature (syntactical and/or semantical), there also exist many pragmatical and methodological reasons for consider combined logics. The major concern in the study of categories of logics (CLE-UNICAMP-Brazil, IST-Lisboa-Portugal) is to describe condition for preservation, under the combination method, of meta-logical  properties (\cite{CCCSS}, \cite{ZSS}). Our complementary approach to this field is study the "global" aspects of categories of logics (\cite{AFLM1}, \cite{AFLM2}, \cite{AFLM3}, \cite{MaMe}).

The initial steps on "global" approach to categories of logics are given in the sequence of papers \cite{AFLM1}, \cite{AFLM2} and \cite{AFLM3}: they present very simple but too strict notions of  logical morphisms, having "good" categorial properties (\cite{AR}) but unsatisfactory treatment of the "identity problem" of logics (\cite{Bez}). More  flexible notions of morphisms between logics are considered in \cite{FC}, \cite{BCC1}, \cite{BCC2}, \cite{CG}: this alternative notion allows better approach to the identity  problem however has many categorial "defects". A "refinement" of those ideas is provided in  \cite{MaMe}: are considered categories of logics satisfying {\em simultaneously} certain natural  conditions: (i) represent the major part of logical systems; (ii) have  good categorial properties; (iii) allow  a natural notion of  algebraizable  logical system (\cite{BP}, \cite{Cze1}); (iv) allow satisfactory treatment of the "identity problem" of logics.
Every category above has the same objects: the  propositional finitary logics, i.e., a pair given by a signature and Tarskian consequence relation on its formula algebra; the morphisms considered are  (some kind of) "logical translations", i.e. some  functions that preserves consequence relations.

 Generalizing the ideas that  describe a precise connection between Boolean algebra and classic propositional logic presented by $Lindenbaum-Tarski$, Blok and Pigozzi  introduced in \cite{BP} the concept of {\em algebraizable logic}, for the first time, as a mathematical definition based on the notions  of algebraizing pair and equivalent algebraic semantics. Here, another relevant  category of logics has, as objects,  the algebraizable logics; the morphisms between them are the translations that preserves algebraizing pairs.

In this work we establish some (categorial) relations between logics and its categories of structures, for instance, given a morphism of algebraizable logics, there is a induced functor between the category of all structures over the underlying signatures such that it restricts to the quasivarieties that are their equivalent algebraic semantics.  About this relation: (i)  we  establish an anti-isomorphism between the class of morphisms of signatures and some functors between the categories of associated structures; (ii) we prove that this anti-isomorphism restricts to  an anti-isomorphism between morphisms of (Lindenbaum) algebraizable logics and some functors on its categories structures that restricts to its quasivarieties.

In section 2 we provide the definitions and basic results (most of them are known, but someone seems new) on categories of logics and algebraizable logics that we will need in the sequence. Our main results are established in section 3 (Theorems \ref{Sf-iso} and \ref{Af-iso}). The last section is devoted to mention possible applications and related works in progress (\cite{MaPi1}, \cite{MaPi2}, \cite{AJMP}).

\section{Preliminaries}

The appearance of several process of combining of logics were the main motivations for the systematic study of categories of logics. Here the objects are logics (signature and consequence operator pairs), the morphisms are translations between logics. In order to find a categories of logics that satisfies some natural requirements, appeared different definitions of categories of logics, more precisely, different definitions of morphisms between logic systems.

\begin{Df}
A signature is a sequence of pairwise disjoint sets $\Sigma=(\Sigma_{n})_{n\in \N}$. In what follows, $X=\{x_{0},x_{1},...,x_{n},...\}$ will denote a fixed enumerable set (written in a fixed order). Denote $F(\Sigma)$ (respectively $F(\Sigma)[n]$), the set of $\Sigma$-formulas over $X$ (respectively, containing exactly the variables $\{x_{0},...,x_{n-1}\}$). Note that the sequence of sets $(F(\Sigma)[n])_{n\in \N}$ is another signature.

A Tarskian consequence relation is a relation $\vdash\subseteq\wp(F(\Sigma))\times F(\Sigma)$, on a signature $\Sigma=(\Sigma_{n})_{n\in \N}$, such that, for every set of formulas $\Gamma,\Delta$ and every formula $\varphi,\psi$ of $F(\Sigma)$, it satisfies the following conditions:
\begin{itemize}
\item[$\circ$]$\bf{Reflexivity:}$If $\varphi\in\Gamma,$ then  $\Gamma\vdash\varphi$;
\item[$\circ$]$\bf{Cut:}$If $\Gamma\vdash\varphi$ and for every $\psi\in\Gamma,\ \Delta\vdash\psi$, then $\Delta\vdash\varphi$;
\item[$\circ$]$\bf{Monotonicity:}$If $\Gamma\subseteq\Delta$ and $\Gamma\vdash\varphi$, then $\Delta\vdash\varphi$;
\item[$\circ$]$\bf{Finitarity:}$If $\Gamma\vdash\varphi$, then there is a finite subset $\Delta$ of $\Gamma$ such that $\Delta\vdash\varphi$;
\item[$\circ$]$\bf{Structurality:}$If $\Gamma\vdash\varphi$ and $\sigma$ is a substitution (see \ref{obs1}), then $\sigma[\Gamma]\vdash\sigma(\varphi)$.
\end{itemize}
\end{Df}

\vtres

The notion of logic that we consider is:

\begin{Df} \label{cons}
A logic of type $\Sigma$, or a $\Sigma-logic$, is a pair $(\Sigma,\vdash)$ where $\Sigma$ is a signature and $\vdash$ is a Tarskian consequence relation.

The set of all consequence relations on a signature $\Sigma$, denoted by $Cons_{\Sigma}$, is endowed with the partial order: $\vdash_0\leq \vdash_1$ iff for each $\Gamma\subseteq F(\Sigma)$, $\{\varphi\in F(\Sigma);\ \Gamma\vdash_{0}\varphi\}=\overline{\Gamma}^{0}\subseteq \overline{\Gamma}^{1}=\{\varphi\in F(\Sigma);\ \Gamma\vdash_{1}\varphi\}$.
\end{Df}

\begin{Obs}
For each signature $\Sigma$, the poset $(Cons_{\Sigma},\leq)$ is a complete lattice. It is in fact an algebraic lattice where the compact elements are the "finitely generated logics", i.e., the logics over $\Sigma$ given by a finite set of axioms and a finite set of (finitary) inference rules.
\end{Obs}

\subsection{Categories of signatures and logics with strict morphisms}

Here we recall the definitions and basic results on the category of signatures with "strict" morphism $\Ss$ and its associated category of logics $\Ls$,  according to \cite{AFLM1}, \cite{AFLM2} and \cite{AFLM3}.

\begin{Df}
The objects of the category $\Ss$ are signatures. If $\Sigma,\Sigma'$ are signatures, then a morphism $f:\Sigma\to\Sigma'$ is a sequence of functions $f=(f_{n})_{n\in \N}$, where $f_{n}:\Sigma_{n}\to\Sigma'_{n}$. For each morphism $f:\Sigma\to\Sigma$ there is only one function $\hat{f}:F(\Sigma)\to F(\Sigma')$, called the extension of $f$, such that:
\begin{itemize}
\item[$\circ$] $\hat{f}(x)=x$, if $x\in X$ ($X$ is a fixed enumerable set);
\item[$\circ$] $\hat{f}(c_{n}(\psi_{0},...,\psi_{n-1}))= f _{n}(c_{n})(\hat{f}(\psi_{0}),...,\hat{f}(\psi_{n-1}))$, if $c_{n}\in \Sigma_{n},\ n\geq 0$
\end{itemize}

Then, by induction, $\hat{f}(\varphi(\psi_{0},...,\psi_{n-1})=\hat{f}(\varphi)(\hat{f}(\psi_{0}),...,\hat{f}(\psi_{n-1}))$, for each $\Sigma$-formula $\varphi$.
\end{Df}

The categories $\Ss$ and $Set^{\N}$ are equivalent, thus we have that $\Ss$ has good categorial properties, namely $\Ss$ is a finitely locally presentable category, in particular it is complete and cocomplete, and the finitely presentable (fp) signatures are the "finite support" signatures (i.e., $\bigcup_{n \in \N} \Sigma_n$ is finite).
\vtres

\begin{Obs}\label{obs1}
\begin{itemize}
\item[(i)] For any substitution function $\sigma:X\to F(\Sigma)$, there is only one extension $\tilde{\sigma}:F(\Sigma)\to F(\Sigma)$ such that $\tilde{\sigma}$ is an homomorphism $\tilde{\sigma}(x) = \sigma(x)$, for all $x\in X$ and $$\tilde{\sigma}(c_{n}(\psi_{0},...,\psi_{n-1})=c_{n}(\tilde{\sigma}(\psi_{0}),...,\tilde{\sigma}(\psi_{n-1}))$$ for all $c_{n}\in \Sigma_{n},n\in\N$. The identity substitution induces the identity homomorphism on the formula algebra; the composition substitution of the substitutions $ \sigma,\sigma':X\to F(\Sigma)$ is the substitution $\sigma'':X\to F(\Sigma),\ \sigma''=\sigma\star\sigma':=\tilde{\sigma}\circ\sigma$ and $\tilde{\sigma''}=\widetilde{\sigma\star\sigma'}=\tilde{\sigma}\circ\tilde{\sigma}'$.
\vtres

Let $f:\Sigma\to \Sigma'$ be a $\Ss$-morphism. Then for any substitution $\sigma:X\to F(\Sigma)$ there is another substitution $\sigma'$ such that $\tilde{\sigma}'\circ \hat{f}=\hat{f}\circ \tilde{\sigma}$.

\[\xymatrix{
F(\Sigma)\ar[rr]^{\hat{f}}\ar[dd]_{\tilde{\sigma}}&&F(\Sigma')\ar[dd]^{\tilde{\sigma}'}\\
&\circlearrowright&\\
F(\Sigma)\ar[rr]_{\hat{f}}&&F(\Sigma')
}\]
\item[(ii)]Let $f:\Sigma\to\Sigma'$ and $\theta\in F(\Sigma)$. If $var(\theta)\subseteq\{x_{i_{0}},...,x_{i_{n-1}}\}$, then $$\hat{f}(\theta(\vec{x})[\vec{x}|\vec{\psi}])=\hat{f}(\theta(\vec{x}))[\vec{x}|\hat{f}(\vec{\psi})].$$ Moreover $var(\hat{f}(\theta))=var(\theta)$ and then $\hat{f}$ restricts to maps $\hat{f}\!\!\rest_{n}:F(\Sigma)[n]\to F(\Sigma')[n], n \in \N$.
\end{itemize}
\end{Obs}
\vtres

Now we introduce the definition of category of logics with "strict" morphism $\Ls$.

\begin{Df}
The objects of $\Ls$ are $l=(\Sigma,\vdash)$, where $\Sigma$ is a signature and $\vdash$ is a tarskian consequence operator. A $\Ls$-morphism, $f:l\to l'$ is a (strict) signature morphism $f\in \Ss(\Sigma,\Sigma')$ such that $\hat{f}:F(\Sigma)\to F(\Sigma')$ is a $(\vdash,\vdash')$-translation: $\Gamma\vdash\psi\Rightarrow\hat{f}(\Gamma)\vdash'\hat{f}(\psi)$
\end{Df}

$\Ls$ is a $\omega$-locally presentable category and the fp logics are given by a finite set of "axioms" and "inference rules" over a fp signature.

Between the categories $\Ls$ and $\Ss$ there exist a forgetful functor $U_s$ such that forgets the consequence relation. $U_s$ has left and right adjoints thus preserves limits and colimits. Moreover $U_s$ "lift" limits and colimits from $\Ss$ to $\Ls$.

\vtres

The categories above mentioned have good categorial properties, but unsatisfactory treatment of the logic problems. e.g., the "identity problem" of logics (\cite{Bez}). Two presentation of classic propositional logic with signatures $\{\neg,\to\}$ and $\{\neg,\vee\}$ do not admits strict morphism between them (because any such morphism must takes $\to$ to $\vee$ and they does not preserve $\vdash$) while it was expected that these presentations would be isomorphic.

\subsection{Categories of signatures and logics with flexible morphisms}

Here we consider the  categories of signatures and of logics based on the developments found in \cite{JKE} \cite{FC}, \cite{BCC1}, \cite{BCC2} and \cite{CG}. This definitions allows a better approach to the "identity problem" of logics, but has many categorial defects.

Similarly to the previous case, firstly we define the category of signature with "flexible" morphism $\Sf$. Before to define this category, let us introduce the following notation:

If $\Sigma=(\Sigma_{n})_{n\in \N}$ is a signature, then $T(\Sigma):=(F(\Sigma)[n])_{n\in \N}$ is a signature too. For each signature $\Sigma$ and $n\in \N$, let the function:

$
\begin{array}{lll}
(j_{\Sigma})_{n}:&\Sigma_{n}\to& F(\Sigma)[n]\\
&c_{n}\mapsto&c_{n}(x_{0},...,x_{n-1})
\end{array}
$

A {\bf flexible} morphism $f : \Sigma\to\Sigma'$  corresponds with a sequence of functions  $f_n^\sharp : \Sigma_n \to F(\Sigma')[n], {n\in\omega}$.

Thus we have the inverse bijections (just notations): \\
$$h\in\Sf(\Sigma,\Sigma')\leftrightsquigarrow h^{\sharp}\in\Ss(\Sigma,T(\Sigma'));
f\in\Ss(\Sigma,T(\Sigma'))\leftrightsquigarrow f^{\flat}\in\Sf(\Sigma,\Sigma').$$

For each morphism $f:\Sigma\to\Sigma'$ in $\Sf$, there is only one function $\check{f}:F(\Sigma)\to F(\Sigma')$, called the extension of $f$, such that:
\begin{itemize}
\item[(i)]$\check{f}(x)=x$, if $x\in X$;
\item[(ii)]$\check{f}(c_{n}(\psi_{0},...,\psi_{n-1}))=f(c_{n})(x_{0},...,x_{n-1})[x_{0}|\check{f}(\psi_{0}),...,x_{n-1}|\check{f}(\psi_{n-1})]$, if $c_{n}\in \Sigma_{n},\\ n \in \N$.
\end{itemize}

\begin{Df}
The category $\Sf$ is the category of signature and $flexible$ morphism as above. The composition in $\Sf$ is given by $(f'\bullet f'')^{\sharp}:=(\check{f}\rest\circ f^{\sharp})$. The identity $id_{\Sigma}$ in $\Sf$ is given by $(id_{\Sigma})^{\sharp}:=((j_{\Sigma})_{n})_{n\in\N}$
\end{Df}

As well as the category $\Ss$, we have that the category $\Sf$ satisfies the conditions of \ref{obs1}.

\begin{Obs} \label{Sf-obs} In \cite{MaMe} is shown that the categories $\Ss$ and $\Sf$ in a way are associated, i.e., there is a pair of adjoint functors between them, namely $(+)_{S}:\Ss\to\Sf$ and $(-)_{S}:\Sf\to\Ss$. Moreover there is a monad or triple $\cT_S = (T_{S}, \mu_{S},\eta_{S})$ on $\Ss$ canonically associated with this adjunction such that $T$ preserve filtered colimits, reflects isomorphisms and, mainly, that $Kleisli(\cT) = \Sf$ (\cite{Mac}). From these results are derived  some additional information on the category $\Sf$: \\
$\bullet$ \ The functor $(+)_S$ preserves colimits and the functor $(-)_S$ preserves limits.\\
$\bullet$ \ Both the functors $(+)_S$ and $(-)_S$ reflect epimorphisms and monomorphisms.\\
$\bullet$ \ The category $\Sf$ has colimits for any (small) diagram "in $\Ss$", i.e., given $\cI$ a small category and a diagram $D : \cI \ \to \ \Ss$, the category $\Sf$ has a colimit for the diagram $(+) \circ D : \cI \ \longrightarrow \ \Sf$. In particular, $\Sf$ has all (small) coproducts and all (small) pushouts "based in $\Ss$".
\end{Obs}

\begin{Df}
The category $\Lf$ is the category of propositional logics and flexible translations as morphisms. This is a category "built above" the category $\Lf$ , that is, there is an obvious forgetful functor $U_{f} :\Lf\to\Sf$.

If $l=(\Sigma,\vdash),l'=(\Sigma',\vdash')$ are logics, then a flexible translation morphism $f:l\to l'$ in $\Lf$ is a flexible signature morphism $f :\Sigma\to\Sigma'$ in $\Sf$ such that "preserves the consequence relation", that is, for all $\Gamma\cup\{\psi\}\subseteq F(\Sigma)$, if $\Gamma\vdash\psi$ then $\check{f}[\Gamma]\vdash'\check{f}(\psi)$. Composition and identities are similar to $\Sf$ .
\end{Df}
\vtres

Due to flexible morphism, this category allows better approach to the identity problem of logics. Consider the flexible morphisms $t:(\rightarrow, \neg)\longrightarrow (\vee',\neg')$ such that $t(\rightarrow)=\neg' x\vee'y$ formula in two variables, $t(\neg)=\neg'$ and $t':(\vee',\neg')\longrightarrow (\rightarrow, \neg)$ such that $t'(\vee')=\neg x\rightarrow y$, $t'(\neg')=\neg$. This pair of morphisms induces an \emph{equipollence} (see \cite{CG})  between these presentations of classic logics. However this category does not has good categorial properties.

\vtres

\begin{Obs}
It follows easily from the facts above that the forgetful functor $U_f:\Lf\to \Sf : ((\Sigma,\vdash)\to((\Sigma ',\vdash '))\mapsto(\Sigma\to\Sigma ')$ has left and right adjoint functors: the left adjoint $\perp_f:\Sf\to \Lf$ and the right adjoint $\top_f:\Sf\to\Lf$ take a signature $\Sigma$ to, respectively, $\perp_f(\Sigma)=(\Sigma,\vdash_{min})$ (the first element of $Cons_{\Sigma}$) and $\top_{f}(\Sigma)=(\Sigma,\vdash_{max})$ (the last element of $Cons_{\Sigma}$). Moreover, $U_{f}\circ\perp_{f}=Id_{\Sf}=U_{f}\circ \top_{f}$ and $U_f$ preserves all limits and colimits that exists in $\Sf$.
\end{Obs}

\begin{Obs}
(a) \ It is known that $\Lf$ has weak products, coproducts and some pushouts, and in the Remark above we see that $U_f$ preserves limits and colimits. As $U_{f}$ also "lift" limits and colimits - the constructions in $\Lf$ are analogous to in $\Ls$  presented in \cite{AFLM3} (just replace everywhere $\hat{f}$ by $\check{f}$) - then given a small category $\mathcal{I}$, $\Lf$ is $\mathcal{I}$-complete (respectively, $\mathcal{I}$-cocomplete) if and only if $\Sf$ is $\mathcal{I}$-complete (respectively, $\mathcal{I}$-cocomplete). Thus, by the Remark \ref{Sf-obs}, it follows  that $\Lf$ has colimits for any (small) diagram "in $\Ls$" (i.e., obtained via $(+) : \Ss \ \to \ \Sf$), in particular, it has all unconstrained fibrings (= coproducts) and the constrained fibrings (= pushouts)  "based in $\Ls$".\\
(b) \ The signature monad $\cT_S= (T_S, \mu_S, \eta_{S})$ associated to the signature adjunction $(\eta_{S},\varepsilon_{S})$ (i.e.,$\mu_{S}=(-)_{S}\varepsilon_{S}(+)_{S}$) "lifts" to a logic monad $\cT_L = (T_L, \mu_{L}, \eta_L)$ associated to the signature adjunction $(\eta_{L},\varepsilon_{L})$ (i.e.,$\mu_{L}=(-)_L\varepsilon_L(+)_L$) and is such that $Kleisli(\cT_L)=\mathcal{L}_{f}$ . Moreover, the functors $(+)_L$ and $(-)_L$ are precisely the canonical functors associated to the adjunction of the Kleisli category of a monad.
\end{Obs}

\subsection{Other categories of logics}

Due to some difficult that was found in the categories of logics mentioned above, are presented in \cite{MaMe}  others categories of logics that  overcome these "defects".

\begin{Fat} \label{othercat}

\begin{enumerate}

\item[(I)] Still on the category $\Lf$ we have the "congruential" logics $\Lf^{c}$. This category is a subcategory of $\Lf$ where the logics are congruential, i.e., logics that satisfies: \[\varphi_{0}\dashv\vdash\psi_{0},...,\varphi_{n-1}\dashv\vdash\psi_{n-1}\Rightarrow c_{n}(\varphi_{0},...,\varphi_{n-1})\dashv\vdash c_{n}(\psi_{0},...,\psi_{n-1}).\]

    The inclusion functor $\Lf^{c}\hookrightarrow\Lf$ has a left adjoint given by congruential closure operator $l \mapsto l^{(c)}$.

A morphism  $f : l \ra l' \in \Lf$ is called {\em dense}, when
$\forall \varphi'_n \in F(\Sigma')[n]$ \ $\exists \varphi_n \in F(\Sigma)[n]$ such that\
 $\varphi'_n \dashv'\vdash \check{f}(\varphi_n)$. If $l' \in \mathcal{L}^{c}_{f}$, then
 $f$ is dense iff $\forall c'_n \in \Sigma'_n$ \ $\exists \varphi_n \in F(\Sigma)[n]$ such that\
  $c'_n(x_0, \ldots, x_{n-1}) \dashv'\vdash \check{f}(\varphi_n)$.

\item[(II)] On the category $\Lf$, consider $Q\Lf$  the quotient category by the congruence relation\footnote{I.e., this category has the same class of objects that $\Lf$, and an arrow between $l \to l'$ the logics is an equivalence class of $\Lf$-arrows $f : l \to l'$.}: $f,g\in \Lf(l,l')$, $f\sim g\ iff\ \check{f}(\varphi)\dashv'\vdash\check{g}(\varphi)$. Thus,  by Proposition 4.3 in \cite{CG}, two logics $l,l'$ are equipollent if only if $l$ and $l'$ are $Q\Lf-$isomorphic. All presentation of classical logic are $Q\Lf-$isomorphic.

\item[(III)] In \cite{MaMe} we found the category $Q\Lf^{c}$ (or simply $\Qf^{c}$).

 For $h \in \mathcal{L}^{c}_{f}(l,l')$,  $[h] \in \Qc(l,l')$ is $\Qc$-isomorphism iff   $h$ is a dense morphism and $h$ is a conservative translation\footnote{I.e., $\Gamma \vdash \psi \ \Leftrightarrow \ \check{h}[\Gamma] \vdash' \check{h}(\psi)$, for all $\Gamma \cup\{\psi\} \subseteq F(\Sigma)$.}.

 This category of logics satisfies {\em simultaneously} certain natural conditions:\\
{\em (a)} \ represent the major part of logical systems; \\
{\em (b)} \ have a good categorial approach (e.g., they are complete, cocomplete and accessible categories); \\
{\em (c)} \ allow a natural notion of algebraizable logical system (\cite{BP},\cite{Cze1}); \\
{\em (d)} \ allow satisfactory treatment of the "identity problem" of logics.\\

\end{enumerate}

\end{Fat}

\subsection{Categories of algebrizable logics}


Traditionally algebraic logic has focused on the algebraic investigation of particular classes of algebras related, in some way, to   logics, whether or not they could be connected to some known assertional system by means of the Lindenbaum-Tarski method. However, when such a connection could be established,
there was interest in investigating the relationship between various meta-logical properties of the logical system and the algebraic properties of the
associated class of algebras.

The  Lindenbaum-Tarski method of algebrization of a logic,  associate a convenient  quotient of  the formula algebra of the logic, by the congruence relation of interprovability: this idea works in classical logic and in some systems of intuitionistic and modal logics.  However this method cannot algebraize other logics. Thus in the end of the 1980's, Blok-Pigozzi  (\cite{BP}) provide a general definition that, in some sense, encompass the traditional method.

Henceforth "algebraizable logic" will mean "algebraizable logic in the Blok-Pigozzi sense".

\begin{Df}
Let $\Sigma$ be a signature. We will denote $\Sigma-Str$ the category with class the objects given by all structures (or algebras) on the signature $\Sigma$ with $\Sigma$-homomorphisms between them. A fundamental example of $\Sigma$-structure is  $F(\Sigma)$, the absolutely free $\Sigma$-algebra on the set $X$.
\end{Df}

\begin{Df}
Given a class of algebras $\textbf{K}$ over the algebraic similarity type $\Sigma$,  the equational consequence associated with $\textbf{K}$ is the relation $\models_{\textbf{K}}$ between a set of equations $\Gamma$ and a single equation $\varphi\equiv \psi$  over $\Sigma$ is defined by:
\[\Gamma\models_{\textbf{K}}\varphi\equiv \psi\ iff\ for\ every\ A\in\textbf{K}\ and\ every\ \Sigma-homomorphism \ \ h:F(\Sigma)\to A,\]
\[\ if \ h(\eta)=h(\nu)\ for\ all\ \eta\equiv\nu\in \Gamma,\ then\ h(\varphi)=h(\psi).\]
\end{Df}

\begin{Df} \label{BP-def}
Let $l=(\Sigma,\vdash)$ be a logic and $\textbf{K}$ be a class of $\Sigma-$algebra. $\textbf{K}$ is a equivalent algebraic semantics to $l$ if $\vdash$ can be faithfully interpreted in $\models_{\textbf{K}}$ of the following sense:
\begin{enumerate}
\item[(1)]there is a finite set $\tau(p)=\{(\delta_{i}(p),\epsilon_{i}(p)),i=1,...,n\}$ of equations in a single variable $p$ such that for all $\Gamma\cup\{\varphi\}\subseteq F(\Sigma)$ and for $j<n$ has been:

    $\Gamma\vdash\varphi\Leftrightarrow\{\tau(\gamma):\gamma\in\Gamma\}\models_{\textbf{K}}\tau(\varphi)$.
\item[(2)]there is a finite system $\Delta_{j}(p,q),j=1,...,m$ of two variables formulas (formed by derived binary connectives) such that for all equation $\varphi\equiv\psi,$\\
$\varphi\equiv\psi=\mid_{\textbf{K}}\models\tau(\varphi\Delta\psi)$\\
where $\varphi\Delta\psi=\Delta(\varphi,\psi)$;  $\Delta(\varphi,\psi)$ abbreviates $\Delta_{j}(\varphi,\psi),j=1,...,m$; $\tau(\varphi\Delta\psi)$ abbreviates $\delta(\varphi\Delta\psi) \equiv \epsilon(\varphi\Delta\psi)$.
\end{enumerate}
In this case we shall say that a logic $l$ is algebraizable. The set $\langle\tau(p),\Delta(p,q)\rangle$ (or just $\langle\tau,\Delta\rangle$) is called an "algebraizing pair", with $\tau = (\delta, \epsilon)$ as the "defining equations" and $\Delta$ as the "equivalence formulas".
\end{Df}

\begin{Fat}\label{BP-def-alter}
Let $K$ an equivalent algebraic semantic for the algebrizable logic $a=(\Sigma,\vdash)$ with algebraizing pair $\langle\tau,\Delta\rangle$, then:
\begin{enumerate}
\item[1.]For all set of equations $\Gamma$ and for all equation $\varphi\equiv\psi$, we have that \[\Gamma\models_{K}\varphi\equiv\psi\ \Leftrightarrow\ \{\xi\Delta\eta:\xi\equiv\eta\in\Gamma\}\vdash\varphi\Delta\psi\]
\item[2.]For each $\psi\in F(\Sigma)$ we have that
\[\psi\dashv \ \vdash\Delta(\tau(\psi)).\]
\end{enumerate}
Conversely, if there is a logic $a=(\Sigma,\vdash)$ and formulas $\langle\Delta(p,q),\tau(p)\rangle$ such that satisfies the conditions 1. and 2., then $K$ is an equivalent algebraic semantics for $a$
\end{Fat}

\begin{Obs} \label{detach} By a direct application of the definition above, if $l=(\Sigma,\vdash)$ is an algebraizable logic and $\phi, \psi \in F(\Sigma)$, then $\phi, \phi \Delta \psi \vdash \psi$ (detachment property).
\end{Obs}

As examples of algebraizable logics we have, in addition to CPC (Classic Propositional Calculus) and IPC (Intuitionistic Propositional Calculus), the modal logics, the Post and Lukasiewicz multi-valued logics, and many of several versions of quantum logic.

In case of CPC (respectively IPC), a possible algebraizing pair $\langle\Delta(p,q),\tau(p) \rangle\\
 =\langle\Delta(p,q), (\epsilon(p),\delta(p))\rangle$ is:

\begin{enumerate}
\item[1.]$\Delta(p,q)=\{p\leftrightarrow q\}$
\item[2.]$\epsilon(p)=p$
\item[3.]$\delta(p)=\top$
\end{enumerate}

and $K$ is the class of Boolean algebras (respectively the class of Heyting algebras). Here the signature of CPC and IPC have as binary connective $\leftrightarrow$. 

Another class of algebras that is the\footnote{See Fact \ref{uniq}.}  equivalent algebraic semantic for an algebrizable logic, but present in many branches of mathematics,  is the class of all groups (\cite{Pa}). To the (equational) theory of groups over the signature $\Sigma=\{\cdot,^{-1},e\}$, it is associated the following propositional logic $l_{Gr}$, the "logic of groups"
\cite{Pa}, over the same signature $\Sigma$, that is

\vsete

Axioms of $l_{Gr}$

\begin{enumerate}
\item[$G_{1}$]$((p\cdot q)\cdot r)\cdot(p\cdot(q\cdot r))^{-1}$
\item[$G_{2}$]$(p\cdot e)\cdot p^{-1}$
\item[$G_{3}$]$(e\cdot p)\cdot p^{-1}$
\item[$G_{4}$]$p\cdot p^{-1}$
\item[$G_{5}$]$p^{-1}\cdot p$
\end{enumerate}

Rules

\begin{enumerate}
\item[$R_{1}$]$p\cdot q^{-1}\vdash q\cdot p^{-1}$
\item[$R_{2}$]$p\cdot q^{-1}\vdash p^{-1}\cdot q^{-1^{-1}}$
\item[$R_{3}$]$\{p\cdot q^{-1},q\cdot r^{-1}\}\vdash p\cdot r^{-1}$
\item[$R_{4}$]$\{p\cdot q^{-1},r\cdot s^{-1}\}\vdash (p\cdot r)\cdot(q\cdot s)^{-1}$
\item[$R_{5}$]$p\vdash p\cdot e^{-1}$
\item[$R_{6}$]$p\cdot e^{-1}\vdash p$
\end{enumerate}

The logic of groups theory has as  algebrizing pair $\langle\Delta(p,q),\tau(p)=\langle\epsilon(p),\delta(p)\rangle\rangle$:
\begin{enumerate}
\item[1.]$\Delta(p,q)=p\cdot q^{-1}$
\item[2.]$\delta(p)=p$
\item[3.]$\epsilon(p)=e$
\end{enumerate}

$K$, in this case, is the class of groups. Worth pointing out that the logic of groups, in some sense, does not admit Deduction Theorem.


\vcinco

Recall that a quasivariety is a class of algebra $K$ such that is axiomatized by quasi-identities, i.e., formulas of the form
\[(p_{1}\equiv q_{1}\wedge...\wedge p_{n}\equiv q_{n})\to p\equiv q\ for\ n\geq 1\]
when $n=0$ the quasi-identity is \[\top\to p\equiv q.\]

Now we will recall a result about "uniqueness" of algebraizing pair and the quasivariety semantics  of an algebraizable logic. For any class $K$ of $\Sigma$-algebras let us denote $(K)^{Q}$ the $\Sigma$-quasivariety generated by $K$.

\begin{Fat} [2.15-\cite{BP}] \label{uniq} Let $a$ be an algebraizable logic.

(a) Let $\langle (\delta_i(p), \varepsilon_i(p)),\Delta_i(p,q)\rangle$, an algebraizing pair for $a$, and $K_{i}$ an equivalent algebraic semantic associated with $a$, for each $i\in \{0,1\}$. Then $(K_{0})^{Q}, (K_{1})^{Q}$ are equivalent algebraic semantics of$a$. Moreover, some uniqueness conditions holds: \\
$\bullet$ \ on quasivariety semantics: $(K_{0})^{Q} = (K_{1})^{Q}$;\\
$\bullet$ \  on equivalence formulas: $\Delta_{0}(p,q)\dashv\vdash\Delta_{1}(p,q)$;\\
 $\bullet$ \ on defining equations: $(\delta_0(p) \equiv  \varepsilon_0(p)) \ =\mid_{K}\models \ (\delta_1(p) \equiv  \varepsilon_1(p))$ (where $K := (K_{0})^{Q} = (K_{1})^{Q}$).

(b) Let $\langle (\delta_i(p), \varepsilon_i(p)),\Delta_i(p,q)\rangle$. Suppose that the following conditions holds:\\
$\bullet$ \  $(\delta_0(p), \varepsilon_0(p)),\Delta_0(p,q)\rangle$ is an algebraizing pair for $a$;\\
$\bullet$ \ $\Delta_{0}(p,q)\dashv\vdash\Delta_{1}(p,q)$;\\
 $\bullet$ \ $(\delta_0(p) \equiv  \varepsilon_0(p)) \ =\mid_{(K_{0})^{Q} }\models \ (\delta_1(p) \equiv  \varepsilon_1(p))$.\\
Then $\langle (\delta_1(p), \varepsilon_1(p)),\Delta_1(p,q)\rangle$is an algebraizing pair for $a$ and $(K_{1})^{Q} = (K_{0})^{Q}$.

\end{Fat}


If $a = (\Sigma, \vdash)$ is an algebraizable logic then, by the Fact above, we can (and we will) denote $QV(a)$ the unique quasivariety on the signature $\Sigma$ that is an equivalent algebraic semantics for $a$.

\begin{Fat} [2.17 \cite{BP}] \label{axiomat}
Let  $a$ be an algebraizable logic $a$ and  $\langle (\delta, \epsilon),\Delta\rangle$ be an algebraizing pair for $a$. Then the quasivariety $QV(a)$ is axiomatized by the set given by the 3 kinds of quasi-equations  below:\\
$\bu$\ $\delta(x_0 \Delta x_0) \equiv \epsilon(x_0 \Delta x_0)$;\\
$\bu$\ $\delta(x_0 \Delta x_1) \equiv \epsilon(x_0 \Delta x_1) \ \to \ x_0 \equiv x_1$;\\
$\bu$\ $(\bigwedge_{i < n} \delta(\psi_i) \equiv \epsilon(\psi_i)) \ \to \  \delta(\phi) \equiv \epsilon(\phi)$, for each $\{\phi, \psi_0, \cdots, \psi_{n-1}\} \sub F(\Sigma)$ such that $\{ \psi_0, \cdots, \psi_{n-1}\} \vdash \phi$, for $n \geq 0$.
\end{Fat}

An attempt to determine if a given logic is algebraizable, at times found difficulties about the definition given above. Thus we have the following characterization.

\begin{Fat} [4.7-\cite{BP}] \label{conditions}
Let $a=(\Sigma,\vdash)$ be a logic and $\Delta\subseteq_{fin} F(\Sigma)[2]$, $(\delta\equiv\epsilon)\subseteq_{fin} (F(\Sigma)[1]\times F(\Sigma)[1])$ such that the conditions below are satisfied
\begin{enumerate}
\item[(a)]$\vdash\varphi\Delta\varphi$, for all $\varphi\in F(\Sigma)$;
\item[(b)]$\varphi\Delta\psi\vdash\psi\Delta\varphi$, for all $\varphi,\psi\in F(\Sigma)$;
\item[(c)]$\varphi\Delta\psi,\psi\Delta\vartheta\vdash\varphi\Delta\vartheta$, for all $\varphi,\psi,\vartheta\in F(\Sigma)$;
\item[(d)]$\varphi_{0}\Delta\psi_{0},...,\varphi_{n-1}\Delta\psi_{n-1}\vdash c_{n}(\varphi_{0},..., \varphi_{n-1})\Delta c_{n}(\psi_{0},...,\psi_{n-1})$, for all $c_{n}\in\Sigma_{n}$ and all $\varphi_{0},\psi_{0},...,\varphi_{n-1},\psi_{n-1}\in F(\Sigma)$;
\item[(e)]$\vartheta\dashv\vdash\Delta(\tau(\vartheta))$, for all $\vartheta\in F(\Sigma)$.
\end{enumerate}
Then $a$ is an algebraizable logic with $\Delta$ as equivalence formulas and $\tau$ as defining equations.

Conversely if $a=(\Sigma,\vdash)$ is a algebrizable logics with algebraizing pair $\langle\Delta(p,q),\\
\tau(p)\rangle$, then the conditions $(a)$ to $(e)$ are satisfied for these formulas.
\end{Fat}

\begin{Obs} \label{strongeralg}
It follows from the characterization above that, if $\vdash_0,  \vdash_1$ are consequence operators over the same signature $\Sigma$, if $l_0 = (\Sigma, \vdash_0)$ is an algebraizable logic with algebraizing pair $\langle\Delta(p,q),\tau(p)\rangle$ and $\vdash_0 \leq  \vdash_1$ (see Definition \ref{cons}), then $l_1 = (\Sigma, \vdash_1)$ is an algebraizable logic and $\langle\Delta(p,q),\tau(p)\rangle$ is an algebraizing pair.
\end{Obs}

With the definition of categories of logics given above, it is possible define categories of algebraizable logics: its morphisms are the translations of algebraizable logics that preserves algebraizing pairs (note that, by Fact \ref{uniq}, this does not depend on particular choice of algebraizing pair of source logic). Other categories of algebraizable logics can be found in \cite{JKE}, \cite{FC}.

\begin{enumerate}
\item[$\bullet$]$\cA_{s}$ is the category of algebraizable logics with morphism in $\Ls$ such that preserves algebraizing pair. In the sequence of works, \cite{AFLM1}, \cite{AFLM2}, \cite{AFLM3} is proven that the category $\cA_s$ is a relatively complete $\omega$-accessible category \cite{AR}.
\item[$\bullet$]$\cA_{f}$ is the category of algebraizable logics with morphisms in $\Lf$ such that preserves algebraizing pair. $\cA_f$ is a (non full) subcategory of $\Lf$, $\cA_{f}\hookrightarrow\Lf$.
\item[$\bullet$] Besides the category $\cA_f$, we consider also the following categories: \\
- \ $\cA_{f}^{c} : = \cA_{f} \cap \mathcal{L}^{c}_{f}$, the (sub)category of algebraizable and congruential logics;\\
- \ $Q\cA_{f}$, the quotient category of $\cA_{f}$ by the congruence determined by interdemonstrability relation ($\dashv \ \vdash$);\\
- \  $Q\cA_{f}^{c}$, the quotient category of $\cA_{f}^c$.
\item[$\bullet$]The "Lindenbaum algebraizable" logics are logics $l\in \cA$ such that given formulas $\varphi,\psi\in F(\Sigma)$, $\varphi\dashv\ \vdash\psi \ \Leftrightarrow \ \vdash\varphi\Delta\psi$ (note this does not depend on the particular choice of $\Delta$; the implication $\Leftarrow$ always hold, by \ref{detach}). The class of Lindenbaum algebraizable logics determines a full subcategory of the category of algebraizable logics ($j:Lind(\cA_{f})\hookrightarrow\cA_{f}$). The category $Lind(\cA_{f})$ is analyzed in the next section of the paper and plays a relevant role in the representation theory of logics (\cite{MaPi1}, \cite{MaPi2}). The inclusion functor $Lind(\cA_{f})\hookrightarrow\cA_{f}$ has a left adjoint functor $L:\cA_{f}\to Lind(\cA_{f})$.
\end{enumerate}

\begin{Df} \label{Deltadense}
(a) Let $l' = (\Sigma', \vdash') \in \Lf$, $a = (\alpha, \vdash) \in \Af$ and $f : l' \to a$ be a $\Lf$-morphism.  Suppose $a \in \cA_f$, then $f$ is called $\Delta-$dense when, given $n \in \N$, for each and  $\varphi\in F(\alpha)[n]$ there is a $\varphi'\in F(\Sigma)[n]$ such that $\vdash\check{f}(\varphi')\Delta\varphi$, for some equivalence formula $\Delta$ of $a$.
Obviously, if $a \in  Lind(\cA_{f})$, then a morphism $f\in \Lf(l',a)$ is $\Delta$-dense iff it is dense.

(b) Let $l = (\Sigma, \vdash) \in \Lf$, $a' = (\alpha', \vdash') \in \Af$. Define the binary relation $\approx_{a'}$ in the set $\Lf(l,a')$
by, $g_0,g_1 \in \Lf(l,a')$,
$$g_0  \approx_{a'} g_1 \ iff \ \forall \phi \in F(\Sigma)(X) \vdash' \check{g}_0(\phi) \Delta' \check{g}_1(\phi),$$
where $\Delta'$ is any equivalence formula for $a'$. It follows from Fact \ref{conditions} that this is an equivalence relation.

When $a = (\alpha', \vdash') \in \Af$, we have an equivalence relation ${}_{a}\!\approx_{a'}$ in the set $\Af(a,a')$. Moreover by the definition of morphisms in $\Af$, the family $\{ {}_{a}\!\approx_{a'} : a, a' \in \Af\}$ defines a congruence relation\footnote{If $f \in Af(b,a),  f' \in \Af(a',b')$, then $(f' \circ g_0 \circ f) {}_{b}\!\approx_{b'} (f' \circ g_1 \circ f)$.} on the category $\Af$ (see \cite{Mac}, Chapter II, Section 8). Denote $\overline{\cA_f}$ the quotient category. It is clear that $\overline{Lind(\cA_f)} = Q(Lind(\cA_f))$.
\end{Df}

By Fact \ref{conditions}, clearly $Lind(\cA_{f}) \subseteq \cA_{f}^{c}$. In the sequence, we establish the equality between these categories \footnotetext{We thank prof. Ramon Jansana for suggesting this result.}. In particular, we obtain that the left adjoint functor $L:\cA_{f}\to Lind(\cA_{f})$ of the inclusion $Lind(\cA_{f})\hookrightarrow\cA_{f}$ is simply given by $l \in \cA_f \mapsto l^{(c)} \in \cA_f^c$ (see Fact \ref{othercat}.(I) and Remark \ref{strongeralg}).

Recall the following

\begin{Df}  Let $\Sigma$ be a signature, ${A}$ be a $\Sigma$-algebra  and $F\subseteq A$.

(a) Let $l = (\Sigma, \vdash)$ be a logic over $\Sigma$. $F$ is called a  $l$-filter if, for all $\Gamma \cup \{\varphi\} \subseteq F(\Sigma)$, if $\Gamma \vdash \phi$ \ $\Rightarrow$ \  $\Gamma \models_A \phi$ (i.e., for each  $\Sigma$-homomorphism $h : F(X) \to A$, if $h[\Gamma] \subseteq F$, then $h(\phi) \in F$.)

(b) Let $\theta$ a congruence in ${A}$. $\theta$ is said to be compatible with $F$ if, for all $a,b\in A$, $a\in F$ and $\langle a,b\rangle\in\theta$ \ $\Rightarrow$ \  $b\in F$.
\end{Df}

\begin{Fat} \label{comp-fa}
(a) [1.5-\cite{BP}]
For any algebra ${A}$ and any $F\subseteq A$, $\Omega_{A}F$ is the largest congruence of ${A}$ compatible with $F$. Where
$\Omega_{A}F = \{\langle a,b\rangle: \varphi^{A}(a,c_{0},...,c_{k-1})\in F \ \Leftrightarrow
 \varphi^{A}(b,c_{0},...,c_{k-1})\in F,$
 for all $\varphi\in Fm_{L}\ and \ c_{i}\in A\}$

(b) [5.2-\cite{BP}]
Let $l = (\Sigma, \vdash)$ be an algebraizable deductive system over the language $\Sigma$, and let $\Delta(x_0,x_1)$ be a system of equivalence formulas. Then \[\Omega_{A}F=\{\langle a,b\rangle:a\Delta^{A}b\in F\}\]
for every $\Sigma-$algebra ${A}$ and every $l-$filter of ${A}$.
\end{Fat}

\begin{Prop}
Let $l$ be a logic. Then $l$ is Lindenbaum algebraizable iff it is an algebraizable and congruential logic.
\end{Prop}

\Dem

``$\Rightarrow$'' Suppose $l \in Lind(\cA_f)$ By Fact \ref{conditions}, it follows that for  every equivalence set of formulas $\Delta$ associated to $l$, the relation defined by $\vdash\Delta(\varphi,\psi)$ is a congruence relation. Therefore that the relation $\dashv \ \vdash$ is a congruence, thus $l \in \cA_f^c$.

``$\Leftarrow$''  Suppose $l \in \cA_f^c$ and let $\varphi,\psi\in F(\Sigma)$. We only have to prove  $\varphi \dashv \ \vdash\psi$ entails $\vdash\varphi\Delta\psi$ (see Remark \ref{detach}).

Consider $T := \{ \gamma  \in F(\Sigma) : \ \vdash \gamma\}$  the set of all theorems of $l$. Let $\varphi,\psi\in F(\Sigma)$ be such that $\varphi\dashv \ \vdash\psi$. Then $\varphi\in T$ iff $\psi\in T$. Thus $\dashv \ \vdash$ is a $\Sigma$-congruence compatible with $T$. Due to the Fact \ref{comp-fa}.(a) above, $\dashv \ \vdash\subseteq\Omega T$,  thus $\langle\varphi,\psi\rangle\in\Omega T$.

It is straightforward that $T$ is a filter in $F(\Sigma)$. By the Fact \ref{comp-fa}.(b) above, $\Omega T=\{\langle \sigma,\sigma'\rangle:\sigma\Delta\sigma'\in T\}$. Therefore $\varphi\Delta\psi\in T$, which means  $\vdash\varphi\Delta\psi$. Therefore $l \in Lind(\cA_f)$.
\qed
\vcinco

We finish this section with the following diagram, that represent the functors (and its adjoints) between some of the categories mentioned above:

\[\xymatrix{
\cA_{f}\ar[rr]^{incl}\ar@<1ex>[dd]^{L}&&\mathcal{L}_{f}\ar[rr]^{q}\ar@<1ex>[dd]^{c}&&\cQ_{f}\ar@<1ex>[dd]^{\bar{c}}\\
&&&&\\
Lind(\cA_{f})\ar[rr]_{incl}\ar@<1ex>[uu]^{j}&&\mathcal{L}_{f}^{c}\ar[rr]_{q^{c}}\ar@<1ex>[uu]^{i}&&\cQ_{f}^{c}\ar@<1ex>[uu]^{\bar{i}}
}
\]

\section{Relations between logics and structures}

This section contains our main contributions. In the tree subsections below, we present: (i)  results on certain adjoint pairs of functors between quasivarieties; (ii) some results about functors between quasivarieties associated to morphisms of (Lindenbaum) algebraizable logics; (iii) a complete (functorial) codification of morphisms of signatures and  of morphisms of algebraizable logics\footnote{In \cite{AJMP}, we provide a encoding of logical morphisms based on different ideas.}.

\vtres

In the sequence: (i) a quasivariety $\cK$ on the signature $\Sigma$ will be viewed as a full subcategory of the category of all structures on that given signature; (ii) for an algebraizable logic $a=(\alpha,\vdash)$, we will denote by $QV(a)$ the unique quasivariety semantics  associated to $a$ (see Fact \ref{uniq}).

\subsection{Quasivarieties and signature functors}

Here we analyze (adjoint pairs of) functors between quasivarieties associated to combination of two fronts: (i) inclusion functors: $\cK \hookr \Sigma-Str$; (ii) "signature" functors i.e. each a $\Sf$-morphism,  $h : \Sigma \lra \Sigma'$, induces a functor $h^\star : \Sigma'-Str \lra \Sigma-Str$. Natural transformations associated to the above mentioned adjunctions also play a significant role here and in the next subsections.

Recall that, by a classical result in universal algebra due to  Mal'cev, a subclass $\cK \sub \Sigma-Str$ is a quasivariety iff it is closed under isomorphisms, substructures, products and ultraproducts (or directed colimits).

\begin{Lem} \label{adjQV-le} Let $\cK$ be a quasivariety on the signature $\alpha$. The inclusion functor has a left adjoint $(L,I): \cK \rightleftarrows \alpha-Str$: given by $M \mapsto M/\theta_M$ where $\theta_M$ is the least $\Sigma$-congruence in $M$ such that $M/\theta_M \in \cK$. Moreover, the unity of the adjunction $(L,I)$ has components $(q_M)_{M \in \Sigma-Str}$, where $q_M : M \thra M/\theta_M$ is the quotient homomorphism.
\end{Lem}

\Dem
Consider $\Gamma_{M}=\{\theta\subseteq|M|\times|M|;$ is congruence relation and $M/\theta \in \cK\}$. $\Gamma$ is not empty, because $\theta=|M|\times|M|$ is a congruence relation and $M/\theta = \{ \star\} \in \cK$. Let $\theta_{M}=\bigcap\Gamma_{M}$.  We will show first that $\theta\in \Gamma_{M}$: as  $\theta_M$ is a $\Sigma$-congruence in $M$, it remains to check that $M/\theta_M \in \cK$.


Consider the "diagonal" $\Sigma$-homomorphism: \[ \delta_M : M\to\prod_{\theta\in\Gamma_{M}}M/\theta;\ m\mapsto([m]_{\theta})_{\theta\in\Gamma_{M}}.\] We will show that $Ker(\delta_M)=\theta_{M}$:

$(m,n)\in Ker(\delta_M)\Leftrightarrow ([m]_{\theta})_{\theta\in \Gamma_{M}}=([n]_{\theta})_{\theta\in \Gamma_{M}}\Leftrightarrow [m]_{\theta}=[n]_{\theta}, \forall\ \theta\in\Gamma_{M}\Leftrightarrow m\theta n\ \forall\ \theta\in\Gamma_{M}\Leftrightarrow m\theta_{M} n.$

Thus, by the "theorem of homomorphism" on $\Sigma-Str$, there is a unique $\Sigma$-\underline{mono}morphism $\bar{\delta}_M : M/\theta_M \rightarrowtail \prod_{\theta\in\Gamma_{M}}M/\theta$ such the diagram below commutes
\[\xymatrix{
M\ar[r]^(.3){\delta_M}\ar[d]_{q_M}&\prod_{\theta\in\Gamma_{M}}M/\theta\\
M/\theta_{M}\ar[ur]_{\bar{\delta}_M}
}\]

As $\cK$ is closed under products, we have that $\prod_{\theta\in\Gamma_{M}}M/\theta \in \cK$. We also have that $\cK$ is closed under substructures and isomorphisms, then $M/\theta_{M}$ is $\cK$.

Denote $L(M) := M/\theta_M$. We will show that $q_M : M \twoheadrightarrow I(L(M))$ satisfies the universal property relatively to $\Sigma$-homomorphisms $f : M \lra I(N)$, with $N \in \cK$.

Thus we obtain a injective $\Sigma$-homomorphism $\bar{f} : M/Ker(f) \rightarrowtail I(N)$. As $\cK$ is closed by substructures and isomorphisms, so we have that $M/Ker(f) \in \cK $. Hence $Ker(f)\in\Gamma_{M}$ and $\theta_M \sub Ker(f)$. Then, again by the theorem of homomorphism, there is a unique homomorphism $\tilde{f} : M/\theta_M \lra N$ such that the following diagram commutes

\[\xymatrix{
M \ar[r]^{q_M}\ar[dr]_{f}&I(L(M))\ar[d]^{I(\tilde{f})}\\
&I(N)
}\]





It follows from  an well known result on adjunct functors, see for instance \cite{Mac}, Theorem 2 in page 81, that there is a unique way to obtain a functor  $L: \Sigma-Str\to \cK$ such that $(q_M)_{M \in \Sigma-Str}$ become the unity of an adjunction $(L,I): \cK \rightleftarrows \alpha-Str$. Given $g \in \Sigma-Str(M, P)$, then $q_P \circ g \in \Sigma-Str(M, I(P))$ and, as $(g\times g)^{-1}[\theta_P] = Ker(q_P \circ g)$, we have that $L(g) = \widetilde{q_P \circ g} : M/\theta_M  \lra P/\theta_P$ : $[m]_{\theta_{M}} \mapsto [g(m)]_{\theta_{P}}$. \qed









\begin{Obs}\label{adjQV-re} Let $\Sigma$ be a signature and  $\cK \sub \Sigma-Str$ be a quasivariety.

(a) The  forgetful functor $(\Sigma-Str \overset{U}\to Set)$ has the "absolutely free algebra" functor $(Set \overset{F}\to \Sigma-Str)$, $Y \mapsto F(Y)$, as left adjunct. The unity of this adjunction has components the inclusion maps $\sigma_Y : Y \rat U(F(Y))$, for each set $Y$.

(b)The (forgetful) functor $(\cK \overset{I}\to \Sigma-Str \overset{U}\to Set)$ has the (free) functor $(Set \overset{F}\to \Sigma-Str \overset{L}\to \cK)$, $Y \mapsto F(Y)/\theta_{F(Y)}$, as left adjunct. Moreover, if $\sigma_Y : Y \ra U\circ F(Y)$ is the $Y$-component of the  unity of the adjunction $(F,U)$, then $(Y \overset{t_Y}\to UILF(Y))  \ : = \ (Y \overset{\sigma_Y}\to UF(Y)  \overset{U(q_{F(Y)})}\to U I L F(Y))$ is the $Y$-component of the adjunction $(L\circ F, U \circ I)$.
\end{Obs}

\begin{Prop}\label{lindalg-ct}
Let $a = (\Sigma, \vdash)$ be an algebraizable and consider the binary relation on $F(X)$,
$$\phi \sim_\Delta \psi \ iff  \ \vdash \phi\Delta\psi,$$
where $\Delta$ is an equivalence formula for $a$. Then:\\
(a) \ $\sim_\Delta$ is a $\Sigma$-congruence on $F(X)$.\\
(b) \ $F(X)/\Delta := F(X)/\sim_\Delta \in QV(a)$.\\
(c) \ $\sim_\Delta = \theta_{F(X)}$ (see Lemma \ref{adjQV-le}), thus $F(X)/\Delta = L(F(X))$ is the free $QV(a)$-object over the set $X = \{x_0, \ldots, x_n, \ldots\}$.\\
In particular, when $a$ is a Lindenbaum algebraizable logic, $F(X)/\Delta  = F(X)/\!(\dashv \ \vdash)$ is the free $QV(a)$-object over the set $X$.

\end{Prop}

\Dem

{\bf (a)} By items (a)-(d) in Fact \ref{conditions} is clear that $\sim_\Delta$ is a $\Sigma$-congruence on $F(X)$.

{\bf (b)} By (a) above,  thus $F(X)/\Delta := F(X)/\sim_\Delta$  is a $\Sigma$-structure.
Thus, to obtain $F(X)/\Delta \in QV(a)$, it is enough to show that $F(\Sigma)/\Delta$ satisfies the conditions of Fact \ref{axiomat}.

$\bu$ \ Let $\varphi := x_0\Delta x_0$, then $\vdash\varphi$. As $a$ is algebraizable logics, $\varphi\dashv\vdash\delta(\varphi)\Delta\varepsilon(\varphi)$. So $\vdash\delta(\varphi)\Delta\varepsilon(\varphi)$. Therefore $[\delta(\varphi)]_{\Delta}=[\varepsilon(\varphi)]_{\Delta}$. Hence $F(X)/\Delta\models \delta(\varphi)\equiv \varepsilon(\varphi)$.
\vtres

$\bu$ \ Suppose $F(X)/\Delta \models \delta(x_0 \Delta x_1) \equiv \epsilon(x_0 \Delta x_1)$. Then
 \[[\delta(x_0 \Delta x_1)]_{\Delta}=[\varepsilon(x_0 \Delta x_1)]_{\Delta}\]
therefore $\vdash \delta(x_0 \Delta x_1) \Delta \epsilon(x_0 \Delta x_1) \to \ x_0 \equiv x_1$. As $a$ is an algebraizable logic, $(x_0 \Delta x_1) \dashv \ \vdash \delta(x_0 \Delta x_1) \Delta \epsilon(x_0 \Delta x_1) \to \ x_0 \equiv x_1$, we obtain $\vdash x_0 \Delta x_1$, i.e. $[x_0]_\Delta = [x_1]_\Delta$. Hence $F(X)/\Delta\models (x_0 \equiv x_1)$ and $F(X)/\Delta \models \delta(x_0 \Delta x_1) \equiv \epsilon(x_0 \Delta x_1) \ \to \ x_0 \equiv x_1$.
\vtres

$\bu$ \  Given $\psi_{0},...,\psi_{n-1},\varphi\in F(X)$ such that $\{\psi_{0},...,\psi_{n-1}\}\vdash\varphi$ and suppose $F(X)/\Delta\models\delta(\psi_{0})\equiv\varepsilon(\psi_{0})\wedge...\wedge\delta(\psi_{n-1})\equiv\varepsilon(\psi_{n-1})$. Then $[\delta(\psi_{0})]_{\Delta}=[\varepsilon(\psi_{0})]_{\Delta},...,[\delta(\psi_{n-1})]_{\Delta}=[\varepsilon(\psi_{n-1})]_{\Delta}$. Therefore
\[\vdash\delta(\psi_{0})\Delta\varepsilon(\psi_{0}),...,\vdash\delta(\psi_{n-1})\Delta\varepsilon(\psi_{n-1}).\]
As $a$ is algebraizable logic, $\psi_i \dashv \ \vdash\psi_{i}$, $ \forall i<n$, thus $\vdash \psi_0,...,\vdash\psi_{n-1}$ and,  by cut, we obtain $\vdash\varphi$. Again, as $a$ is algebraizable, we obtain $\vdash \delta(\varphi)\Delta \varepsilon(\varphi)$. Hence $F(X)/\Delta\models \delta(\varphi)\equiv \varepsilon(\varphi)$ and $F(X)/\Delta\models (\bigwedge_{i < n} \delta(\psi_{i})\equiv\varepsilon(\psi_{i}))  \ \to \ \delta(\varphi)\equiv \varepsilon(\varphi)$.

{\bf (c)} Let $M \in QV(a)$. The universal property of $\sigma_X : X \lra U(F(X))$ induces a bijection $\Sigma-Str(F(X), I(M)) \cong Set(X, U(I(M))$: for each function $v : X \lra U(I(M))$ there is an unique $\Sigma$-homomorphism $V : F(X) \lra I(M)$ such that $V \circ \sigma_X = v$. Establish the equality $\sim_\Delta = \theta_{F(X)}$ is equivalent to  prove that $\sim_\Delta \sub Ker(V)$, for each function  $v : X \lra U(I(M))$. Suppose $\phi \sim_{\Delta} \psi$,  then $\vdash \phi \Delta \psi$. As $a$ is an algebraizable logic we obtain, by Fact \ref{BP-def-alter} $\models_{QV(a)} \phi \equiv \psi$, i.e. for each $M \in QV(a)$ and each $\Sigma$-homomorphism $H : F(X) \lra I(M)$, $H(\phi) = H(\psi)$. Thus $\sim_{\Delta} \sub Ker(V)$  for each function  $v : X \lra U(I(M))$.
\qed

\begin{Obs} \label{freerest} By   reasoning  analogous to in proof above we can establish that, for every $Y\subseteq X$, the binary relation on $F(Y)$ given by $(\sim_\Delta)\rest := (\sim_\Delta) \cap (F(Y)\times F(Y))$ coincides with $\theta_{F(Y)}$, thus $F(Y)/\Delta\rest := F(Y)/\sim_{\Delta}\rest$ is the free $QV(a)$-object over the set $Y$.

\end{Obs}

\begin{Ct} \label{h*} {\bf Signature functors:} Given a morphism in $\Sf$,  $\Sigma \overset{h}\to \Sigma'$, we associate a functor $\Sigma-Str \overset{h^\star}\leftarrow \Sigma'-Str$ in the following way

$\bu$ \ For each $M' \in \Sigma'-Str$ denote $h^{\star}(M') = (M')^h$ the $\Sigma$-structure such that\\
-- $|(M')^h| = |M'|$ (structures with same underlying set);\\
-- Let $k \geq 0$ and $c_k \in (\Sigma)_k$, then $h(c_k) \in F(\Sigma')[k]$  is a {\em first-order $k$-ary  term over $\Sigma'$} and its  interpretation in the $\Sigma'$-structure $M'$ is a certain $k$-ary operation on $|M'|$,   ${M'}^{{h(c_k)}} : |M'|^k \to |M'|$;  define  $(c_k)^{(M')^h} := h(c_k)^{M'}$  (it is a $k$-ary operation on $|M'^h|$).

If $\phi \in F(\Sigma)$ has exactly $n$ variables, then it can be viewed as n-ary first-order $\Sigma$-term and its interpretation over $(M')^h$ is defined (by recursion on complexity); analogously  the n-ary first-order $\Sigma'$-term $\check{h}(\phi)$ can be interpreted on $M'$.   We can prove, by induction on the complexity of $\phi$, that the $n$-operations on the same set $|(M')^h| = |M'|$,  $(\phi)^{(M'^{h})}, {(\check{h}(\phi))}^{(M')}$, coincide.

$\bu$ \ Let $g\in \Sigma-Str(M',N')$, we define $h^{\star}(M',g,N')= (M'^{h}, g , N'^h) \in \Sigma-Str(M'^h,N'^h)$: clearly, the function $g$ determines a $\Sigma$-homomorphism from $M'^h$ into $N'^h$).

It is clear that $h^\star$ preserves identities and composition, thus it is a (covariant) functor.

By construction, the functor $h^\star : \Sigma'-Str \rightarrow \Sigma-Str$   "commutes over $Set$", i.e., $U \circ h^\star = U'$. It is straitforward that $h^\star$ preserves, {\em strictly}, the following constructions: substructures, products, directed inductive limits,  reduced products, congruences and quotients.


\end{Ct}

\begin{Prop} \label{adjunct-prop} Consider a signature morphism $h \in \cS_f(\Sigma,\Sigma')$ and quasivarieties $I : \cK \hookrightarrow \Sigma-Str$, $I' : \cK' \hookrightarrow \Sigma'-Str$. Suppose that the induced functor $h^\star : \Sigma'-Str \lra \Sigma-Str$  restricts  to a $h^\star\!\!\!\rest : \cK' \rightarrow \cK$, i.e. there is a (unique) functor  $h^\star\!\!\!\rest$ such that $I \circ h^\star\!\!\!\rest = h^\star \circ I'$, then

(a) $h^\star\!\!\!\rest : \cK' \rightarrow \cK$ has a left adjunct $G : \cK \to \cK'$.

(b) Suppose that $h^\star\!\!\!\rest : \cK' \rightarrow \cK$ satisfies the following conditions\footnote{A first kind of  examples is given  by inclusion functors $I : \cK \hookr \Sigma-Str$, for some quasivariety $\cK \sub \Sigma-Str$. We will see more examples in section 3.2 below and produce a classification of all "logical examples" in section 3.3. In \cite{MaPi2}, we will apply these results.}:\\
(b1) $h^\star\!\!\!\rest$ is faithful;\\
(b2) $h^\star\!\!\!\rest$ is full;\\
(b3) $h^\star\!\!\!\rest$ is injective on objects;\\
(b4) $h^\star\!\!\!\rest$  is  {\em hereditary}, i.e., given $M\in \cK$, $N' \in \cK'$ such that there is an injective $\Sigma$-homomorphism $j : M \rat h^{\star}\!\!\rest(N')$,  then there is $M' \in \cK'$ such that $h^{\star}\!\!\rest(M') = M$.\\
Then the left adjunct $G$ can be defined on objects $M \in \cK$ as "a quotient" $G(M) \in \cK'$, with   $h^\star\rest(G(M)) = M/\rho_M$, where $\rho_M$   is the least $\Sigma$-congruence in $M$ such that $M/\rho_{M} = h^\star(M')$, for some $M' \in \cK'$ (that is automatically unique by (l3)); moreover the $M$-component of the unity of the adjunction is the quotient map $p_M : M \thra M/\rho_M$.

\end{Prop}

\Dem (a) We will give here an indirect proof of the existence of the left adjunct $G$: we will prove that the hypothesis on "Freyd Left Adjoint  Theorem" (see \cite{Mac}, Theorem 2, page 117) are satisfied by $h^\star\!\!\!\rest : \cK' \rightarrow \cK$.

$\bu$ \ As $\cK \sub \Sigma-Str$ and $\cK' \sub \Sigma'-Str$ are closed under isomorphisms, substructures and products, $\cK$ and $\cK'$ are complete categories, i.e. they have all small limits. Moreover, as $h^\star : \Sigma'-Str \to \Sigma-Str$ (strictly) preserves: isomorphisms, substructures  and products, then the same holds for $h^\star\!\!\rest : \cK' \rightarrow \cK$. Thus $h^\star\!\!\rest : \cK' \rightarrow \cK$ preserves all small limits.

$\bu$ \ We show that the "solution set condition" holds for $h^\star\!\!\rest$. Let $M \in \cK$ and consider $\kappa := card(|M|)$ and consider the {\em class} $C_M := \{ N' \in \cK':$  such that $N'$ has a $\cK'$-generator subset of size $\leq \kappa \}$. It is clear that there is a {\em set} $S_M \sub C_M$ of representatives of $C_M$ modulo isomorphism. We will show that $\bigcup_{S' \in S_M} \cK(M, h^\star\!\!\rest (S'))$ is a {\em set} that satisfies the solution set condition for $M'$.\\
Let $P' \in \cK'$ and $f : M \to h^\star\!\!\rest(P')$ be a $\cK$-morphisms. Let $N' \sub P'$ be the $\Sigma'$-substructure of $P'$ that is generated by $image(f)$. Then $N' \in C_M$ and we can take $S' \in S_M$ such that $S' \cong_{\cK'} N'$. Consider a fixed $\cK'$-isomorphism $t : S' \to N'$ and let $i : N' \hookr P'$ be the inclusion. Then we have shown that the homomorphism $f : M \to  h^\star\!\!\rest(P')$ factors through some member $g$ of the set $\bigcup_{S' \in S_M} \cK(M, h^\star\!\!\rest (S'))$ (i.e. $f = h^{\star}(i \circ t) \circ g$).

(b) Let $M \in \cK$, we will prove that the "quotient map" $p_M : M \to h^\star\!\!\rest(G(M))$, $h^\star\!\!\rest(G(M)) = M/\rho_M$, satisfies the universal property.
Consider $\Omega_M := \{ \theta \sub M \times M: \theta$ is a $\Sigma$-congruence in $M$ and there is a (unique) $P'_\theta \in \cK'$ such that $M/\theta = h^\star\!\!\rest(P'_\theta)\}$. We show first that $\Omega_M$ has  minimum by verifying that $\rho_M := \bigcap \Omega_M    \in \Omega_M$. Indeed we have a {\em injective} $\Sigma$-homomorphism $j: M/ \rho_M \rat \prod_{\theta \in \Omega_M} M/\theta$, $[m]_{\rho_M} \mapsto ([m]_{\theta})_{\theta \in \Omega_M}$. By definition of $\Omega_M$, $\prod_{\theta \in \Omega_M} M/\theta = \prod_{\theta \in \Omega_M} h^\star\!\!\rest(P'_\theta)$.  As $h^\star\!\!\rest$ preserves products we have the {\em injective} $\Sigma$-homomorphism $j : M/\theta_M \rat h^\star\!\!\rest(\prod_{\theta \in \Omega_M} P'_\theta)$. By conditions (b4) and (b3), $M/\rho_M = h^\star\!\!\rest(M')$ for a unique $M' \in \cK$. Thus $\rho_M = \bigcap \Omega_M \in \Omega_M$.

Let $N' \in \cK'$ and $f : M \to h^\star\!\!\rest(N')$ be a $\Sigma$-homomorphism:  we will show that there is a unique $\Sigma'$-homomorphism $f' : M' \to N'$ such that:
$$ (M \overset{f}\lra h^\star\!\!\rest(N')) \ = \ (M \overset{p_M}\thra h^\star\!\!\rest(M') \overset{h^\star\!\!\rest(f')}\lra h^\star\!\!\rest(N'))$$

   Then $f$ factors through the quotient homomorphism $q_f : M \thra M/Ker(f)$   by the {\em injective} $\Sigma$-homomorphism $\bar{f} : M/Ker(f) \rat h^{\star}\!\!\rest(N')$. Then, by conditions (b4) and (b3), $M/ker(f) = h^{\star}\!\!\rest(P')$ for a unique $P' \in \cK$. As $Ker(f) \in \Omega_M$, we have $\rho_M \sub Ker(f)$ and, by the theorem of homomorphism, there is a unique $\Sigma$-homomorphism $\bar{f} : M/\rho_M \to h^\star\!\!\rest(N')$ such that $\bar{f}\circ p_M = f$. As $M/\rho_M = h^\star\!\!\rest(M')$ for a unique $M' \in \cK$,  the conditions (b1) and (b2) ensures that there is a unique $\Sigma'$-homomorphism $f' : M' \to N'$ such that $h^\star\!\!\rest(f') = \bar{f}$. Then $f'$ is the unique $\Sigma'$-homomorphism such that $f =  h^{\star}\!\!\rest(f') \circ p_M$.
\qed

\begin{Prop}\label{diag-adj} Consider a signature morphism $h \in \cS_f(\Sigma,\Sigma')$ and quasivarieties $I : \cK \hookrightarrow \Sigma-Str$, $I' : \cK' \hookrightarrow \Sigma'-Str$. Suppose that the induced functor $h^\star : \Sigma'-Str \rightarrow \Sigma-Str$  restricts  to a (unique) functor $h^\star\!\!\rest : \cK' \rightarrow \cK$, i.e. $I \circ h^\star\!\!\rest = h^\star \circ I'$. Denote $G$ and $\bar{G}$ the (unique up to natural isomorphism) left adjunct functors of, respectively , $h^\star$ and $h^\star\rest$ (they exists by Proposition \ref{adjunct-prop} above). Then:

(a)   $(G \circ F) \cong F'$ and $(\bar{G} \circ L)  \cong (L' \circ G)$.

(b) There is a natural epimorphism $\tilde{h}:  L \circ h^\star \thra h^\star\!\!\rest \circ L$, that restricts to $ L \circ h^\star\circ I' = h^\star\!\!\!\rest \circ L' \circ I'$.

\end{Prop}

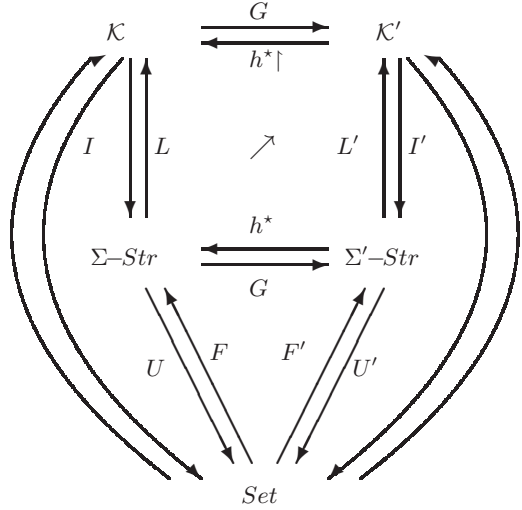
\[\begin{picture}(160,120)
\setlength{\unitlength}{.6\unitlength} \thicklines
\put(0,20){\small $\Sigma\!\!-\!\!Str$}
\put(10,162){\small $\cK$}
\put(180,162){\small $\cK'$}
\put(40,90){\small $L$}
\put(-5,90){\small $I$}
\put(100,145){\small $h^\star\!\!\rest$}
\put(100,175){\small $\bar{G}$}
\put(160,20){\small $\Sigma'\!\!-\!\!Str$}
\put(100,90){$\nearrow$}
\put(155,90){\small $L'$}
\put(200,90){\small $I'$}
\put(100, 40){\small $h^\star$}
\put(100, 0){\small $G$}
\put(25,150){\vector(0,-1){100}}
\put(195,150){\vector(0,-1){100}}
\put(35,50){\vector(0,1){100}}
\put(185,50){\vector(0,1){100}}
\put(150,160){\vector(-1,0){80}}
\put(150,30){\vector(-1,0){80}}
\put(70,170){\vector(1,0){80}}
\put(70,20){\vector(1,0){80}}
\put(35,-50){\small $U$}
\put(165,-50){\small $U'$}
\put(75,-40){\small $F$}
\put(120,-40){\small $F'$}
\put(95,-130){\small $Set$}
\put(102,-105){\vector(-1,2){55}}
\put(118,-105){\vector(1,2){55}}
\put(35,05){\vector(1,-2){55}}
\put(185,05){\vector(-1,-2){55}}
\qbezier(20,150)(-100,15)(70,-115)
\put(65,-110){\vector(1,-1){05}}
\qbezier(220,150)(340,15)(170,-115)
\put(155,-110){\vector(-1,-1){05}}
\qbezier(0,150)(-120,15)(50,-115)
\put(05,150){\vector(2,1){05}}
\qbezier(200,150)(320,15)(150,-115)
\put(215,150){\vector(-2,1){05}}
\end{picture}\]
\vspace{2cm}

\Dem

(a)  The uniqueness up to isomorphism of left adjuncts entails that $ U \circ h^{\star}$ has a left adjunct isomorphic to $G \circ F$. As $U \circ h^{\star} = U'$ and $F'$ is a left adjunct of $U'$, again the uniqueness of left adjuncts up to isomorphism ensures that $(G \circ F) \cong F'$. Analogously, from the equality $I \circ h^\star\!\!\!\rest = h^\star \circ I'$, we obtain the natural isomorphism $(\bar{G} \circ L)  \cong (L' \circ G)$.

(b)
Let $M' \in \Sigma'-Str$ and consider the canonical arrow in $\Sigma'-Str$ $q'_{M'} : M' \thra M'/\theta'_{M'} = I'(L'(M'))$. Applying $h^\star$, we obtain the  (surjective) $\Sigma$-homomorphism $h^\star(q'_{M'}) : h^\star(M') \thra h^\star(M'/\theta'_{M'})$ and the induced $\Sigma$-isomorphism
$$\overline{h^\star(q'_{M'})} : h^\star(M')/ker(h^\star(q'_{M'}) \overset{\cong}\lra h^\star(M'/\theta'_{M'}).$$
 As the functor $h^\star$ commutes over $Set$ and (strictly) preserves  substructures and products, then $h^\star(\theta'_{M'})$ is a $\Sigma$ congruence over $h^\star(M')$ and
 $h^\star(M'/\theta'_{M'}) = h^\star(M')/h^\star(\theta'_{M'})$; thus $ker(h^\star(q'_{M'})) = h^\star(\theta'_{M'})$ and  \[\overline{h^\star(q_{M'})} = Id : h^\star(M')/h^\star(\theta'_{M'}) \lra h^\star(M'/\theta'_{M'}).\]
 In particular,  $h^\star(M')/h^\star(\theta'_{M'}) = h^\star(M'/\theta'_{M'}) \in \cK$ and $\theta_{h^\star(M')} \subseteq h^\star(\theta'_{M'})$. Therefore, there is a {\em canonical surjective} $\Sigma$-homomorphism
$$\tilde{h}_{M'} : h^\star(M')/\theta_{h^\star({M'})} \thra h^\star(M')/h^\star(\theta'_{M'}) = h^\star(M'/\theta'_{M'}) :$$
 this defines a $\cK$-morphism $\tilde{h}_{M'} : L(h^\star(M')) \thra h^\star\rest(L'(M'))$.

When $M' \in \cK'$, then $h^\star(M') \in \cK$, $\theta'_{M'} = \Delta_{|M'|}$ and $\theta_{h^\star(M')} = \Delta_{|h^\star(M')|} = h^\star(\Delta_{|M'|}) = h^\star(\theta'_{M'})$. Thus, in this case,  $\tilde{h}_{M'} = Id : L(h^\star(M')) \lra h^\star\rest(L'(M'))$.

If $f' : M' \lra N'$ is a $\Sigma'$-homomorphism, then ${h^\star}(f') : {h^\star}(M') \lra {h^\star}(N')$ is a $\Sigma$-homomorphism. To show  that the diagram below commutes

\[\xymatrix{
M'\ar[dd]_{f'} & L\circ h^{\star}(M')\ar[dd]_{L \circ h^\star (f')}\ar[rr]^{\tilde{h}_{M'}}&
&h_{\rest}^{\star}\circ L'(M')\ar[dd]^{h_{\rest}^{\star}\circ L'(f')}\\
&&&\\
N' & L\circ h^{\star}(N')\ar[rr]_{\tilde{h}_{N'}}& &h_{\rest}^{\star}\circ L'(N')
}\]

it is enough to realize that

$$ h_{\rest}^{\star}(L'(f')) \circ \tilde{h}_{M'} \circ q_{h^\star}(M') = \tilde{h}_{N'} \circ L (h^\star (f')) \circ q_{h^\star}(M'),$$

where $q_{h^\star}(M') : {h^\star}(M') \thra {h^\star}(M')/\theta_{{h^\star}(M')} $ is the canonical {\em surjective} $\Sigma$- homomorphism, but this follows immediately from a diagram chase. Thus $\tilde{h} := (\tilde{h}(M')_{M' \in \Sigma'-Str}$ is a natural transformation.
\qed

\subsection{Algebraizable logics and functors}

In this part of the work, we verify that the general results on the functors between quasivarieties presented in the previous subsection can be applied to functors induced by logical morphisms between algebraizable logics. Are established the first connections between properties of the logical morphisms and the properties of its induced functors.

\begin{Prop}\label{restrict} Let $a = (\alpha, \vdash)$ and $a' = (\alpha', \vdash')$ be algebraizable logics and let $h \in \cA_f(a,a')$. Then the induced functor $h^\star : \alpha'-Str \rightarrow \alpha-Str$   restricts to $h^\star\!\!\!\rest : QV(a') \rightarrow QV(a)$ (i.e. $I \circ h^\star\!\!\!\rest = h^\star \circ I'$).

\end{Prop}


\Dem As $QV(a) \sub \alpha-Str$ and $QV(a') \sub \alpha'-Str$ are full subcategories, it is enough to show that: for each  $M' \in QV(a')$ we have $h^\star(M') \in QV(a)$.

 It follows from  the description of a set of quasi-identities that determines the unique equivalent quasivariety semantics associated to algebraizable logic in Fact \ref{axiomat} it follows that, if $(\Delta, (\delta,\epsilon))$ is an algebraizable pair for $a = (\alpha,\vdash)$, then the set
of quasi-identities $S_a = S_a^0 \cup S_a^1 \cup S_a^2$ axiomatizes $QV(a)$, where:\\
$S_a^0 = \{ \delta(x_0 \Delta x_0) \equiv \epsilon(x_0 \Delta x_0)\}$;\\
$S_a^1 = \{ \delta(x_0 \Delta x_1)\equiv\epsilon(x_0 \Delta x_1)) \ \ra \ x_0 \equiv x_1 \}$;\\
$S_a^2 = \{ (\delta(\psi_0)\equiv\epsilon(\psi_0) \wedge...\wedge \delta(\psi_{n-1})\equiv\epsilon(\psi_{n-1})) \ \ra \ \delta(\varphi)\equiv\epsilon(\varphi)\ :\{\psi_{0},..., \psi_{n-1}\} \vdash \varphi \}$.

 Denote $\mathfrak{h}$ the extension of $h$ to first-order formulas, instead $\check{h}$ that is the extension of $h$ for propositional $\alpha$- formulas (= first-order terms). For instance, $\mathfrak{h}((\delta(\psi_0)\equiv\epsilon(\psi_0) \wedge \ldots \wedge \delta(\psi_{n-1})\equiv\epsilon(\psi_{n-1})) \ \ra \ \delta(\varphi)\equiv\epsilon(\varphi)) = (\check{h}\delta(\check{h}\psi_0)\equiv\check{h}\epsilon(\check{h}\psi_0) \wedge \ldots \wedge \check{h}\delta(\check{h}\psi_{n-1})\equiv\check{h}\epsilon(\check{h}\psi_{n-1})) \ \ra \ \check{h}\delta(\check{h}\varphi)\equiv\check{h}\epsilon(\check{h}\varphi)$.

As $h \in \cA_f(a,a')$ then:\\
$\bu$ \  $((\check{h}(\delta),\check{h}(\epsilon)), \check{h}(\Delta))$ is an algebraizable pair for $a'$. \\
$\bu$ \  If $\{\psi_{0},...,\psi_{n-1}\} \vdash \varphi$, then $\{\check{h}\psi_{0},..., \check{h}\psi_{n-1}\} \vdash' \check{h}\varphi$.\\
From these, it follows that: $\mathfrak{h}[S_{a}^0] = S_{a'}^0$, $\mathfrak{h}[S_{a}^1] = S_{a'}^1$ and  $\mathfrak{h}[S_{a}^2] \sub S_{a'}^2$.
Thus, for each quasi-equation $\Omega \in S_{a}^0 \cup  S_{a}^1 \cup S_{a}^2$, we have $M' \vDash_{\alpha'} \mathfrak{h}(\Omega)$. On the other hand, for each first-order formula $\Theta$ holds the following equivalence:
$$M' \vDash_{\alpha'} \mathfrak{h}(\Theta) \ \Leftrightarrow \ h^\star(M')  \vDash_{\alpha} \Theta.$$
Thus $h^\star(M')  \in QV(a)$, as we wish.
\qed

\begin{Prop} \label{QV-quo-prop} Let  $l = (\Sigma, \vdash) \in \Lf$ and $a, a' \in \Af$. Keeping the notation in the definition \ref{Deltadense}, we have:

(a) \ Let  $g_0, g_1 : l \to a'$ be $\Lf$-morphisms.  Then
$$g_0 \approx_{a'} g_1 \ \Leftrightarrow \ g_0^\star\!\!\rest = g_1^\star\!\!\rest : QV(a') \ra \Sigma-Str.$$

(b) \ Let $g_0, g_1 : a \to a'$ be $\cA_{f}$-morphisms. Then
$$[g_0]_{\approx} = [g_1]_{\approx} \in \overline{\cA_f} \ \Leftrightarrow \ g_0^\star\!\!\rest = g_1^\star\!\!\rest : QV(a') \ra QV(a).$$

\end{Prop}

\Dem Item (b) follows from item (a), since a quasivariety on signature $\alpha$ determines a full subcategory of $\alpha-Str$.

(a)$"\Rightarrow"$

Let $M'\in QV(a')$ and $c_{n}\in \Sigma_{n}$. As $g_0 \approx_{a'} g_1$, we have that
$$\vdash_{a'} \check{g}_{0}(c_{n})(x_{0},...,x_{n-1})\Delta \check{g}_{1}(c_{n})(x_{0},...,x_{n-1}).$$
Thus, by Fact \ref{BP-def-alter},
$\models_{QV(a')}\check{g}_{0}(c_{n})(x_{0},...,x_{n-1})\equiv \check{g}_{1}(c_{n})(x_{0},...,x_{n-1})$.
Therefore:
\[c_{n}^{{M'}^{g_{0}}} = (g_{0}(c_{n}))^{M'} = (g_{1}(c_{n}))^{M'} = c_{n}^{{M'}^{g_{1}}}\]

Thus $g_{0\rest}^{\star}(M') = g_{1\rest}^{\star}(M')$ and, as $g_{0\rest}^{\star}, g_{1\rest}^{\star}$ commute over $Set$, they coincide also on the arrow level. Therefore $g_{0\rest}^{\star} = g_{1\rest}^{\star}$.
\vtres

$"\Leftarrow"$

Suppose that $g_{0\rest}^{\star}=g_{1\rest}^{\star}$. Let $\varphi \in F(\Sigma)$,  hence $\varphi^{{M'}^{g_{0}}}= \varphi^{{M'}^{g_{1}}}$ for all $M'\in QV(a')$.  So $\models_{QV(a')}\check{g}_{0}(\varphi)\equiv \check{g}_{1}(\varphi)$. Due to $a'$ to be algebraizable, by Fact \ref{BP-def-alter}, $\vdash_{a'} \check{g}_{0}(\varphi) \Delta \check{g}_{1}(\varphi)$. Therefore $g_0 \approx_{a'} g_1$.
\qed

\vtres

\begin{Cor} \label{QV-quo-cor} Let  $l = (\Sigma, \vdash) \in \Lf$ and $a, a' \in$ {$\cA^c_f$}.

(a) \ Let  $g_0, g_1 : l \to a'$ be $\Lf$-morphisms.  Then
$$[g_0]_{\dashv \vdash} = [g_1]_{\dashv \vdash} \in \Qf \ \Leftrightarrow \ g_0^\star\!\!\rest = g_1^\star\!\!\rest : QV(a') \ra \Sigma-Str.$$

(b) \ Let $g_0, g_1 : a \to a'$ be $\cA^c_{f}$-morphisms. Then
$$[g_0]_{\dashv \vdash} = [g_1]_{\dashv \vdash} \in Q\cA^c_f  \ \Leftrightarrow \ g_0^\star\!\!\rest = g_1^\star\!\!\rest : QV(a') \ra QV(a).$$

\end{Cor}

\begin{Prop} \label{QV-iso-prop} Let $a$ and $a'$ be  algebraizable logics and $a \overset{h'}{\underset{h}\rightleftarrows} a'$   be a pair of $\cA_{f}$-morphisms. Then $a \overset{[h']_\approx}{\underset{[h]_\approx}\rightleftarrows} a'$ is a pair of inverse $\overline{\cA_f}$-isomorphisms
\ iff \ $QV(a) \overset{{h'}^\star\!\!\rest}{\underset{{h}^\star\!\!\rest}\leftrightarrows} QV(a')$ is an isomorphism of categories.

\end{Prop}

\Dem
 The induced $\overline{\cA_f}$-morphisms $a \overset{[h']_\approx}{\underset{[h]_\approx}\rightleftarrows} a'$ is a pair of inverse $\overline{\cA_f}$-isomorphisms iff
$$[id_{a}]_\approx=[h']_\approx\circ[h]_\approx = [h'\bullet h]_\approx \ and$$
$$ [id_{a'}]_\approx =[h]_\approx \circ[h']_\approx =[h\bullet h']_\approx$$
$$\ iff \ (by \ Corollary\ \ref{QV-quo-prop}.(b))$$
$$id_{a}^\star\!\!\rest= (h'\bullet h)^\star\!\!\rest =  (h^\star \circ  h'^\star)\!\!\rest =  h^\star\!\!\rest \circ  h'^\star\!\!\rest \ and \ id_{a'}^\star\!\!\rest =   (h\bullet h')^\star\!\!\rest =  (h'^\star \circ  h^\star)\!\!\rest = h'^\star\!\!\rest \circ  h^\star\!\!\rest $$
$$iff$$
the pair of functors $QV(a) \overset{{h'}^\star\!\!\rest}{\underset{{h}^\star\!\!\rest}\leftrightarrows} QV(a')$ is a pair of inverse isomorphism of categories.
\qed

Restricting the above result to the setting of Lindenbaum algebarizable logics, we obtain the

\begin{Cor} \label{QV-iso-cor} Let $a$ and $a'$ be Lindenbaum algebraizable logics and $a \overset{h'}{\underset{h}\rightleftarrows} a'$   be a pair of $\cA^c_{f}$-morphisms. Then $a \overset{[h']_{\dashv \vdash}}{\underset{[h]_{\dashv \vdash}}\rightleftarrows} a'$ is a pair of inverse $Q(\cA^c_f)$-isomorphisms\footnote{Remember that $Q(\cA^c_f) = Q(Lind(\cA_f)) = \overline{Lind(\cA_f)}$.}  \ iff \ $QV(a) \overset{{h'}^\star\!\!\rest}{\underset{{h}^\star\!\!\rest}\leftrightarrows} QV(a')$ is an isomorphism of categories.
\end{Cor}

\begin{Prop}\label{QV-dense-prop} Let  $l = (\Sigma, \vdash) \in \Lf$ and $a, a' \in \Af$.

(a) Let  $h : l \ra a'$ be a $\Lf$-morphism. Consider the conditions:

(a1) \  $h$ is a $\Delta$-dense $\Lf$-morphism.

(a2) \ The functor $h^\star\!\!\rest : QV(a') \ra \Sigma-Str$ is full, faithful, injective on objects and satisfies the heredity  condition (see \ref{adjunct-prop}.(b4)).

(b) Let  $h: a \ra a'$ be a $\Af$-morphism. Consider the conditions:

(b1) \  $h$ is a $\Delta$-dense $\Af$-morphism.

(b2) \ The functor $h^\star\!\!\rest : QV(a') \ra QV(a)$ is full, faithful, injective on objects and satisfies the heredity  condition.

Then (a1) $\Rightarrow$ (a2) and (b1) $\Rightarrow$ (b2)
\end{Prop}

\Dem The implication [(b1) $\Rightarrow$ (b2)] follows from  [(a1) $\Rightarrow$ (a2)] and  Proposition \ref{restrict}, since the inclusion functor $I : QV(a) \hookr \alpha-Str$ is clearly full, faithful, injective on objects and satisfies the heredity  condition.

 We will prove [(a1) $\Rightarrow$ (a2)]

 $\underline{Full}$: Let $M', N'\in QV(a')$ and $f : h_{\rest}^{\star}(M')\to h_{\rest}^{\star}(N')$ be a $\Sigma$-homomorphism. As $h_{\rest}^{\star}$ commutes over $Set$, we have $U(f): |M'| \to |N'|$ is a function. We will prove that $f: M'\to N'$ is a $\alpha'$-homomorphism.

By the hypothesis ($h$ is $\Delta$-dense), for each  $c'_{n}\in \alpha'_{n}$  there is $\varphi_{n}\in F(\Sigma)[n]$ such that $\vdash'\check{h}(\varphi_{n}(x_{0},...,x_{n-1})) \Delta' c'_{n}(x_{0},...,x_{n-1})$. Thus,  as $a'$ is an algebraizable logic, then
$$\models_{QV(a)} \check{h}(\varphi_{n}(x_{0},...,x_{n-1}))\equiv c'_{n}(x_{0},...,x_{n-1}).$$

Let $v: X \to |M'|$ be a function. Consider $m_{0}= v(x_{0}),...,m_{n-1} = v(x_{n-1})$. So $\check{h}(\varphi_{n}(x_{0},...,x_{n-1})))^{M'}[\overrightarrow{x}/\overrightarrow{m}]=(c'_n(x_{0},...,x_{n-1}))^{M'}[\overrightarrow{x}/\overrightarrow{m}]$.

$\begin{array}{rcl}
f((c'_n(x_{0},...,x_{n-1}))^{M'}[\overrightarrow{x}/\overrightarrow{m}])&=&f(\check{h}(\varphi_{n}(x_{0},...,x_{n-1})))^{M'}[\overrightarrow{x}/\overrightarrow{m}])\\
&=&f((\varphi_{n}(x_{0},...,x_{n-1}))^{{M'}^{h}}[\overrightarrow{x}/\overrightarrow{m}])\\
&=&(\varphi_{n}(x_{0},...,x_{n-1}))^{{N'}^{h}}[\overrightarrow{x}/f(\overrightarrow{m})]\\
&=&(c'_n(x_{0},...,x_{n-1}))^{{N'}}[\overrightarrow{x}/f(\overrightarrow{m})]
\end{array}$

Therefore $f$ is a $QV(a)$-morphism.

$\underline{Faithful}$: Let $f_1, f_2 \in QV(a')(M',N')$. As $h_{\rest}^{\star}(M',f, N') = ({M'}^h, f, {N'}^h)$, if $h_{\rest}^{\star}(f_{1})=h_{\rest}^{\star}(f_{2}) \in \Sigma-Str({M'}^h,{N'}^h)$ then $f_{1} = f_{2}$.

$\underline{Injective\ on\ objects}$: Let $M', N'\in QV(a')$ such that $h_{\rest}^{\star}(M') = h_{\rest}^{\star}(N')$, so $|M'|=|N'|$. Given $c'_{n}\in\alpha_{n}$ there is $\varphi_{n}(x_{0},...,x_{n-1})\in F(\Sigma)[n]$  such that
$$\vdash'\check{h}(\varphi_{n}(x_{0},...,x_{n-1}))\Delta' c'_{n}(x_{0},...,x_{n-1})\ \Rightarrow$$
$$\models_{QV(a)}\check{h}(\varphi_{n}(x_{0},...,x_{n-1}))\equiv c'_{n}(x_{0},...,x_{n-1})$$

Hence, given $m_{0},...,m_{n-1} \in |M'|$

\[\begin{array}{rcl}
c'^{M'}_{n}(m_{0},...,m_{n-1})&=&\check{h}(\varphi_{n})^{M'}(m_{0},...,m_{n-1})\\
&=&\varphi_{n}^{M'^{h}}(m_{0},...,m_{n-1})\\
&=&\varphi_{n}^{N'^{h}}(m_{0},...,m_{n-1})\\
&=&\check{h}(\varphi_{n})^{N'}(m_{0},...,m_{n-1})\\
&=&c'^{N'}_{n}(m_{0},...,m_{n-1})
\end{array}\]

Therefore $M' = N'$

$\underline{Heredity}$:

Let $M\in \Sigma-Str$, $N' \in QV(a')$ be such that there is an injective $\Sigma$-homomorphism $j : M \rat h^{\star}\!\!\rest(N')$.  We must show that  there is $M' \in QV(a')$ such that $h^{\star}\!\!\rest(M') = M$. Remark that: \\
(i) as $h^\star\rest$ is injective on objects, then $M'$ is unique; \\
(ii) as $h^\star\rest$ is full and faithful, then $j : M' \rat N'$ is an injective $\alpha'$-homomorphism.\\
 Thus, as $N' \in QV(a')$, it is enough construct an $\alpha'$-structure $M'$ such that $h^{\star}(M') = M$, because then $M' \in QV(a')$ automatically.

 As $h$ is $\Delta$-dense, given $c'_{n}\in\alpha_{n}$ {\em select a formula} $\langle {c}'_{n} \rangle(x_{0},...,x_{n-1})\in F(\Sigma)[n]$  such that
$$\vdash'\check{h}(\langle {c}'_{n}\rangle (x_{0},...,x_{n-1}))\Delta' c'_{n}(x_{0},...,x_{n-1})\ \Rightarrow$$
$$ \models_{QV(a')}\check{h}(\langle c'_{n} \rangle(x_{0},...,x_{n-1}))\equiv c'_{n}(x_{0},...,x_{n-1})$$

Define $|M'| := |M|$. Let $m_{0},...,m_{n-1} \in |M'|$, define $c'^{M'}_n(m_{0},...,m_{n-1}) : =$ \\
$\langle c'_{n} \rangle^{M}(m_{0},...,m_{n-1})$. Then:
$$j(c'^{M'}_{n}(m_{0},...,m_{n-1})) = j(\langle c'_{n} \rangle^{M}(m_{0},...,m_{n-1})) =$$
$$ = \langle c'_{n} \rangle^{N'^h}(j(m_{0}),...,j(m_{n-1})) = (\check{h}(\langle c'_{n} \rangle))^{N'}(j(m_{0}),...,j(m_{n-1})) =$$
$$ = c'^{N'}_{n}(j(m_{0}),...,j(m_{n-1})). $$
Thus $j : M' \rat N'$ is an injective $\alpha'$-homomorphism. In particular, $M' \in QV(a')$.

Now let $c_k \in \Sigma_k$ and $a_0, \cdots, a_{k-1} \in |M'|$. Then
 $$j(c_k^{{M'}^h}(a_0, \cdots, a_{k-1})) = j((\check{h}(c_k))^{M'}(a_0, \cdots, a_{k-1})) =$$
$$= (\check{h}(c_k))^{N'}(j(a_0), \cdots, j(a_{k-1})) =c_k^{{N'}^h}(j(a_0), \cdots, j(a_{k-1})) =$$
$$= j(c_k^{M}(a_0, \cdots, a_{k-1})).$$

As $j$ is injective, then $c_k^{{M'}^h}(a_0, \cdots, a_{k-1}) = c_k^{M}(a_0, \cdots, a_{k-1})$.
Hence $h^{\star}\!\!\rest(M') = M$.
\qed

\vtres

As density and $\Delta$-density of morphisms coincide on Lindenbaum algebraizable logics, we immediately obtain the

\begin{Cor}\label{QV-dense-cor} Let  $l = (\Sigma, \vdash) \in \Lf$ and $a, a' \in \cA^c_f$.

(a) Let  $h : l \ra a'$ be a $\Lf$-morphism. Consider the conditions:

(a1) \  $h$ is a dense $\Lf$-morphism.

(a2) \ The functor $h^\star\!\!\rest : QV(a') \ra \Sigma-Str$ is full, faithful, injective on objects and satisfies the heredity  condition .

(b) Let  $h: a \ra a'$ be a $\cA_f^c$-morphism. Consider the conditions:

(b1) \  $h$ is a dense $\cA_f^c$-morphism.

(b2) \ The functor $h^\star\!\!\rest : QV(a') \ra QV(a)$ is full, faithful, injective on objects and satisfies the heredity  condition.

Then (a1) $\Rightarrow$ (a2) and (b1) $\Rightarrow$ (b2)
\end{Cor}

In the  next subsection we will be able to prove that the implications presented in Proposition \ref{QV-dense-prop} and Corollary \ref{QV-dense-cor} are, in fact, equivalences.

\begin{Prop} \label{adj-rest}
 Let $a\overset{h}\to a'\in \cA^c_{f}$, then:

(a)  $h^{\star}_{\rest} : QV(a')\to QV(a)$ has a left adjunct $G : QV(a) \to QV(a')$.

(b) In case that $h$ is a dense morphism, then the left adjunct $G$ can be defined on objects $M \in QV(a)$ as "a quotient" $G(M) \in QV(a')$, with   $h^\star\rest(G(M)) = M/\rho_M$, where $\rho_M$   is the least $\Sigma$-congruence in $M$ such that $M/\rho_{M} = h^\star(M')$, for some (and unique) $M' \in \cK'$ (that is automatically unique by (l3)); moreover the $M$-component of the unity of the adjunction,  $M  \to h^\star(G(M))$, is the quotient map $p_M : M \thra M/\rho_M$.
\end{Prop}

\Dem Item (a) follows from Propositions \ref{restrict} and \ref{adjunct-prop}.(a). Item (b) is a direct consequence of Proposition \ref{adjunct-prop}.(b) and Corollary \ref{QV-dense-cor}.(b).
\qed

\begin{Obs}\label{godel}
Let $a=IPC$ and $a'=CPC$ both Lindenbaum algebraizable logics with the same signature. We have the inclusion morphism $h : IPC\to CPC$. So $h^{\star}_{\rest}:BA\to HA$ has left a adjoint functor $G:HA\to BA$. Observe that $h^{\star}_{\rest}$ is the inclusion functor. Hence given $H\in HA$, $G(H)=H/F_{H}$, where $F_{H}$ is the filter in $H$ generated by the subset $\{a\leftrightarrow \neg\neg a:a\in H\}$.
Its possible to proof that $G(H) \cong H_{\neg\neg}$, where $H_{\neg\neg}$ denote the poset of regular elements of $H$, that is, those elements $x\in H$ such that $\neg\neg x=x$.

This fact motivate us to investigate the relation of  G\"odel translation with the left adjunct functor $G$.

As an application of some of the general results in the present work, we derive in \cite{MaPi2} a generalized "Glinvenko's Theorem" related to an $\cA^c_f$-morphism $h : a \to a'$, whenever holds a simple condition of the unity of the adjunction $(G,h^{\star}_{\rest}) : QV(a') \leftrightarrows QV(a)$.

\end{Obs}

\subsection{Functorial encoding of logical morphisms}

Here we apply the previous results to determine a faithful  representation of algebraizable logic morphisms as certain  functors. We start presenting a  functorial encoding of signature morphisms.

\begin{Lem}\label{sigfunt-le}
 Let $\Sigma, \Sigma' \in Obj(\Sf)$. Consider $H: \Sigma'-Str \ra \Sigma-Str$ a functor that "commutes over $Set$" (i.e. $U \circ H = U'$) and, for each set $Y$, let $\eta_H(Y) : F(Y) \ra H(F'(Y))$ be the unique $\Sigma$-morphism such that $(Y \overset{\sigma_Y}\ra UF(Y)  \overset{U(\eta_Y)}\to UHF'(Y)) \ = \ (Y \overset{\sigma'_Y}\to U'F'(Y))$ (by the universal property of $\sigma_Y$). Then:\\
(a)\ $(\eta_H(Y))_{Y \in Set}$ is a natural transformation $\eta_{H}: F \ra H\circ F'$.\\
(b)\ For each set $Y$ and each $\psi \in F(Y)$, $Var(\eta_H(Y)(\psi)) \ \subseteq \ Var(\psi)$.\\
When $\forall \psi \in F(Y)$, $Var(\eta_H(Y)(\psi)) \ = \ Var(\psi)$, we will say that $\eta_H(Y)$ "preserves variables".\\
(c)\ For each  $n \in \N$, let $X_n := \{x_0, \cdots, x_{n-1}\} \subseteq X$, if $\eta_H(X_n)$ preserves variables, then the mapping $c_n \in \Sigma_n \ \mapsto \ \eta_H(X_n)(c_n(x_0, \cdots, x_{n-1})) \in F'(X_n)$ determines a unique $\Sf$-morphism $m_H : \Sigma \ra \Sigma'$ such that $\check{m}_H = \eta_H(X)$.
\end{Lem}

\Dem

(a)\ \ Let a function $f : Y\to Z$. Consider the diagram:

\[\xymatrix{
Y\ar[r]^{\sigma_{Y}}\ar@/^2pc/[rr]^{\sigma'_{Y}}\ar[d]_{f}&UF(Y)\ar[r]^{U(\eta_{H}(Y))}\ar[d]_{UF(f)}&UHF'(Y)\ar[d]^{UHF'(f)}\\
Z\ar[r]_{\sigma_{Z}}\ar@/_2pc/[rr]_{\sigma'_{Z}}&UF(Z)\ar[r]_{U(\eta_{H}(Z))}&UHF'(Z)
}\]
\vcinco

The left square commutes because $\sigma$ is the unit of adjunction between $U$ and $F$ and the external diagram commutes because $\sigma'$ is the unit of adjunction between $U' =U \circ H$ and $F'$. We also have that $U(\eta_{H}(Y))\circ\sigma_{Y}=\sigma'_{Y}$; the same is valid when we change $Y$ for $Z$. Thus, a diagram chase entails ensures that
$$UHF'(f) \circ U(\eta_H(Y)) \circ \sigma_Y = U(\eta_H(Z)) \circ UF(f) \circ \sigma_Y.$$
As $U$ is a functor, the universal property of $\sigma_Y$ give us
$$ U(HF'(f) \circ \eta_H(Y)) = U (\eta_Z \circ F(f)).$$
As $U$ is faithful, we obatin
 $$ HF'(f) \circ \eta_H(Y) = \eta_Z \circ F(f).$$
 Thus $\eta_{H} : F \Rightarrow H \circ F'$ is a natural transformation.
\vtres

(b)\ \ Let $Y$ a set and $\psi\in F(Y)$. Consider $Z = Var(\psi)$. and denote $i : Z \hookr Y$ the inclusion function.    As $\eta_{H}$ is a natural transformation, we have the follow commutative diagram:

\[\xymatrix{
F(Y)\ar[r]^{\eta_{H}(Y)}&HF'(Y)\\
F(Z)\ar[u]^{F(i)=incl}\ar[r]_{\eta_{H}(Z)}&HF'(Z)\ar[u]_{HF'(i)=incl}
}\]

As $\psi \in F(Z)$,  we have $\eta_{H}(Z)(\psi)=\eta_{H}(Y)(\psi)$. Therefore $Var(\eta_{H}(Y)(\psi))\subseteq Z = Var(\psi)$.
\vtres

(c) Follows directly from the definition of flexible morphism of signatures.
\qed

Note that any functor $H: \Sigma'-Str \ra \Sigma-Str$  that "commutes over $Set$" ($U \circ H = U'$) is automatically faithful, since $U'$ is a faithfull functor. Another simple but  useful result is given by the following

\begin{Fat} \label{preservevar} Keeping the notation above, are equivalent:

(a) $\eta_H(Y)$ preserves variables, for each set $Y$.

(b) $\eta_H(Y)$ preserves variables, for each set $Y \sub X$.

(c) $\eta_H(X_n)$ preserves variables, for each $n  \in \N$.

\end{Fat}

\Dem We only have to show that $(c) \Rightarrow (a)$.

Let $Y$ be an arbitrary set and let $\psi \in F(Y)$. Let  $\{y_0, \cdots, y_{n-1}\}  \sub Y$ be a (bijective) enumeration of  $Var(\psi)$ and consider the injection $x_i \overset{f}\mapsto y_i$, $i < n$, from $X_n$ into $Y$. As $F(f)$ is injective, denote $\tilde{\psi} \in F(X_n)$ the unique member such that $F(f)(\tilde{\psi}) = \psi \in F(Y)$. By hypothesis $ Var(\tilde{\psi}) =  X_n = Var(\eta_H(X_n)(\tilde{\psi}))$. As $F(f)$ and $H(F'(f))$ are injective and $\eta_H$ is a natural transformation, then a diagram chase entails
$$Var(\eta_H(Y)({\psi})) = Var({\psi}).$$
\qed

When $\eta_H$ satisfies the equivalent conditions above, we say that  $\eta_H$ "preserves variables". As we will see in the sequence, this is a fundamental concept in this work, leading us to the following

\begin{Df} \label{sigfunct} Let $\Sigma, \Sigma' \in Obj(\Sf)$ and $H :
\Sigma'-Str \ra \Sigma-Str$ be functor. We will say that $H$ is a
"signature" functor if it satisfies the conditions below:\\
$(s1)$ $H$ commutes over $Set$ (i.e. $U \circ H = U'$);\\
$(s2)$ $\eta_H$ preserves variables.
\end{Df}

\begin{Prop} \label{compsigfun}

(a) Let $\Sigma-Str \overset{id}\to \Sigma-Str$. Then $\eta_{id_{\Sigma}-Str} = id_{F}$ and $id_{\Sigma-Str}$ is a
signature functor;   moreover $m_{id_{\Sigma-Str}} = id_{\Sigma} \in \Sf(\Sigma,\Sigma)$.

(b) Let $(\Sigma-Str \overset{H}\leftarrow \Sigma'-Str \overset{H'}\leftarrow
\Sigma''-Str)$ be functors that commutes over $Set$. Then $\eta_{H\circ H'} = H(\eta_{H'})\circ\eta_{H}$. If $H$ and $H'$ are signature functors, then $H \circ
H'$ is a signature functor and, moreover, in this case,  $m_{H \circ H'} = m_{H'} \bullet m_H \in \Sf(\Sigma, \Sigma'')$.
\end{Prop}

\Dem
(a) It is clear that $id_{\Sigma-Str}$ commutes over $Set$. For each set $Y$, notice that the function $id_{F(Y)} : F(Y) \to id_{\Sigma-Str}(F(Y))$ satisfies $\sigma_{Y}\circ U(id_{F(Y)}) = \sigma_{Y}$ where $\sigma_{Y}: Y \to UF(Y)$ is the unit of the adjunction $F\dashv U$. Then, the universal property of $\sigma_Y$ entails $\eta_{id_{\Sigma-Str}}({Y}) = id_{F(Y)}$. Thus, in particular, $\eta_{id_{\Sigma-Str}}$ preserves variables, i.e., $id_{\Sigma-Str}$ is a signature functor. Moreover, $\check{m}_{id_{\Sigma-Str}} = \eta_{id_{\Sigma-Str}}(X) = id_{F(X)} : F(X) \to F(X)$, thus $m_{id_{\Sigma-Str}} = id_{\Sigma} \in \Sf(\Sigma, \Sigma)$.

\vtres

(b) As $H'$ and $H$ commute over $Set$, we have that $H'\circ H$ also commutes over $Set$. We have that the following commutative diagrams:

\[\xymatrix{
Y\ar[r]^{\sigma_{Y}}\ar[dr]_{\sigma'_{Y}}&UF(Y)\ar[d]^{U(\eta_{H}(Y))}&&&&Y\ar[r]^{\sigma'_{Y}}\ar[dr]_{\sigma''_{Y}}&U'F'(Y)\ar[d]^{U'(\eta_H'({Y}))}\\
&U'F'(Y)&&&&&U''F''(Y)
}\]

As $U \circ H = U'$, we obtain
$$  U (\eta_{H \circ H'}) \circ \sigma_Y = \sigma''_Y  =$$
$$= UH(\eta_{H'}(Y)) \circ U(\eta_H(Y)) \circ \sigma_Y = U( H(\eta_{H'}(Y)) \circ \eta_H(Y) ) \circ \sigma_Y.$$
By the universal property of $\sigma_Y$, for each set $Y$, we obtain
$\eta_{H\circ H'} = H(\eta_{H'})\circ\eta_{H}$.

For each $n \in \N$, $\eta_{H\circ H'}(X_n) = H(\eta_{H'}(X_n)) \circ \eta_H(X_n)$.
Now suppose that $H$ and $H'$ are signature functors. As $\eta_H(X_n)$ and $\eta_{H'}(X_n)$ preserve variables and $H$ commutes over $Set$, then $\eta_{H\circ H'}(X_n)$ preserves variables. Thus $H \circ H'$ is a signature functor. Moreover, in this case, $\check{m}_{H \circ H'} = \eta_{H\circ H'}(X) = H(\eta_{H'}(X)) \circ \eta_H(X) = H(\check{m}_{H'}) \circ \check{m}_H$: this means that $m_{H \circ H'} = m_{H'} \bullet m_{H}$.
\qed

\vcinco

In the sequence, we will see that, among the functors  $H : \Sigma'-Str \to \Sigma-Str$ that commutes over $Set$,  there are two kinds  of functors that also preserves variables: the isomorphisms  $\Sigma'-Str \to \Sigma-Str$  and the functors $h^\star$, induced by $\Sf$-morphisms $h : \Sigma \to \Sigma'$.

\begin{Prop} \label{isosig} Let $H : \Sigma'-Str \to \Sigma-Str$ be an isomorphism of categories such that $U \circ H = U'$. Then $H$ is a signature functor.
\end{Prop}

\Dem As  $H : \Sigma'-Str \to \Sigma-Str$ is an isomorphism of categories such that $U \circ H = U'$, then $H^{-1} :  \Sigma-Str \to \Sigma'-Str$ is an isomorphism of categories and obviously $U' \circ H^{-1} = U$. Let $Y$ be a set and consider $\psi \in F(Y)$. By the Lemma \ref{sigfunt-le}.(b), $Var(\eta_H(Y)(\psi)) \ \sub \ Var(\psi)$). On the other hand, by the Proposition \ref{compsigfun} $id_{F(Y)} = \eta_{H\circ H^{-1}}(Y) = H(\eta_{H^{-1}}(Y))\circ\eta_{H}(Y)$, thus $Var(\psi) = Var(H(\eta_{H^{-1}}(Y)) (\eta_{H}(Y)(\psi))) \sub Var(\eta_{H}(Y)(\psi))$, since we can apply Lemma \ref{sigfunt-le}.(b) to $H^{-1}$ and $H$ commutes over $Set$ .
\qed
\vtres

\begin{Prop} \label{h*pres} Let $h \in \Sf(\Sigma, \Sigma')$, then for all $Y \sub X$,  $\eta_{h^\star}(Y) = \check{h}_{\rest Y} : F(Y) \to F'(Y)^h$. In particular, $\eta_{h^\star}$ preserves variables and $h^\star$ is a signature functor according the Definition \ref{sigfunct}.
\end{Prop}

\Dem Firstly observe that the function $\check{h}: U(F(X)) \to U'(F'(X))$ (see subsection 2.2) is such that $Var(\check{h}(\phi)) = Var(\phi)$, for each $\phi \in U(F(X))$. Thus, for each $Y \sub X$, it restricts to $\check{h}_{\rest Y} : U(F(Y)) \to U'(F'(Y))$ and for each $\varphi \in U(F(Y))$,   $Var(\check{h}_{\rest Y} (\varphi)) = Var(\varphi)$.

Now, remark that $\check{h}_{\rest Y}$ determines a $\Sigma$-homomorphism.
\[\begin{array}{rcl}
\check{h}_{\rest Y}:&F(Y)&\to h^\star(F'(Y))\\
&\varphi&\mapsto\check{h}(\varphi)
\end{array}\]

Clearly, the diagram below commutes:

\[\xymatrix{
Y\ar[r]^{\sigma_{Y}}\ar[dr]_{\sigma'_{Y}}&UF(Y)\ar[d]^{U(\check{h}_{\rest Y})}\\
&Uh^{\star}F'(Y)
}\]

Due to the universal property of $\sigma_Y$, we have $\eta_{h^{\star}}(Y)=\check{h}_{\rest Y}$, for each $Y \sub X$. Thus, by Fact \ref{preservevar}, $\eta_{h^{\star}}$ preserves variables and, as $h^\star$ commutes over $Set$, then $h^\star$ is a signature functor according the Definition \ref{sigfunct}.
\qed

The family of functors $h^\star$, induced by $\Sf$-morphisms $h$, have a nice categorial behavior:

\begin{Prop} \label{compsigmor}

(a) Let $\Sigma \overset{id_{\Sigma}}\to \Sigma$ be the identity $\Sf$-morphism on the signature $\Sigma$. Then  $id_{\Sigma}^\star = id_{\Sigma-Str} \in Cat(\Sigma-Str,\Sigma-Str)$.

(b) Let $(\Sigma \overset{h}\to \Sigma' \overset{h'}\to
\Sigma'')$ be $\Sf$-morphisms. Then $(h' \bullet h)^\star = h^\star \circ h'^\star \in Cat(\Sigma''-Str, \Sigma-Str)$.
\end{Prop}

\Dem Since the functors induced by signature morphisms are faithful and commute over $Set$, we only have to verify the equalities of functors in (a) and (b) at level of the objects.

It is clear that, for each $M \in Obj(\Sigma-Str)$, $M =  id_{\Sigma-Str}(M) = id_{\Sigma}^\star(M)$, establishing item (a).

To prove item (b), note first that, for each $M'' \in Obj(\Sigma''-Str)$,
$$U((h' \bullet h)^\star(M'')) = U''(M'') = U'(h'^\star(M'')) =$$
$$= (U\circ h^\star)(h'^\star(M'')) = U ((h^\star \circ h'^\star) (M'')),$$
Thus, the $\Sigma$-structures $(h' \bullet h)^\star(M'')$ and $(h^\star \circ h'^\star) (M'')$ shares the same underlying set.
It remains verify that, for each $n \in \N$ and each $c_{n} \in \Sigma_n$,
$$(I): \ (c_n)^{(h' \bullet h)^\star(M'')} = (c_n)^{(h^\star \circ h'^\star) (M'')}.$$
Developing the left hand side of $(I)$ we obtain
$$(L): \ (c_n)^{(h' \bullet h)^\star(M'')} = ((h' \bullet h) (c_n))^{M''} = ((h' \bullet h) (c_n))^{M''} =$$
$$ = ((\check{h}'\circ h)(c_n))^{M''} = ((\check{h}'( h(c_n)))^{M''} = (h(c_n))^{M''^{h'}}.$$

Developing the right hand side of $(I)$ we obtain
$$(R): \  (c_n)^{(h^\star \circ h'^\star) (M'')} =  (c_n)^{(h^\star (h'^\star (M''))} =$$
$$=   (c_n)^{(h^\star ((M'')^{h'})} = (h(c_n))^{M''^{h'}}.$$

Summing up, we obtain $(h' \bullet h)^\star(M'') = (h^\star \circ h'^\star) (M'')$. Thus $(h' \bullet h)^\star = (h^\star \circ h'^\star)$.
\qed







At this point, is natural consider the following

\begin{Df} \label{Sf+} Let $\Sf^\dagger$ denote the (non-full) subcategory of the category of all the (large) categories\footnote{I.e., the category whose objects are large categories and the
arrows are functors between categories.}  given by  the categories
$\Sigma-Str$, for each signature $\Sigma$,  and with the signature functors as morphisms between them.
\end{Df}

\begin{Teo} \label{Sf-iso}
 The categories $\Sf$ and $\Sf^\dagger$ are anti-isomorphic. More precisely:\\
(a) \ The mapping $\Sigma \in Obj(\Sf) \mapsto \Sigma-Str \in Obj(\Sf^\dagger)$ is
bijective; \\
(b) \ Given $\Sigma, \Sigma' \in \Sf$, the mappings $h \in
\Sf(\Sigma,\Sigma') \ \mapsto\ h^\star \in
\Sf^\dagger(\Sigma'-Str,\Sigma-Str)$ and $H \in
\Sf^\dagger(\Sigma'-Str,\Sigma-Str)\ \mapsto\ m_H \in
\Sf(\Sigma,\Sigma')$   are (well defined) inverse bijections.\\
(c) \ $id_{\Sigma}^\star = id_{\Sigma-Str}$ and $(h' \bullet h)^\star = h^\star \circ h'\star$;\\
$m_{id_{\Sigma-Str}} = id_{\Sigma}$ and $m_{H \circ H'} = m_{H'} \bullet m_H$.
\end{Teo}

\Dem The equalities in item (c) were established in Propositions \ref{compsigfun} and \ref{compsigmor}.

(a) The mapping $\Sigma \in Obj(\Sf) \mapsto \Sigma-Str \in Obj(\Sf^\dagger)$ is surjective, by definitions of  $\Sf^\dagger$.
 Note that $\Sigma \neq \Sigma' \ \Rightarrow \ \Sigma-Str \neq \Sigma'-Str$ (in fact, $\Sigma \neq \Sigma' \ \Rightarrow \ Obj(\Sigma-Str) \cap Obj(\Sigma'-Str) = \emptyset$).

(b) The mappings $H \mapsto m_H$ and $h \mapsto h^\star$   are well defined by, respectively, Lemma \ref{sigfunt-le}.(c) and Proposition \ref{h*pres}. Moreover, these results ensures that $\check{m}_{h^\star} = \eta_{h^\star}(X) = \check{h}$. Therefore $m_{h^\star} = h$.
It  remains only to prove that, for each signature functor $H : \Sigma'-Str \to \Sigma-Str$, $(m_H)^\star = H$.

It is enough to prove that $H(M') = (m_H)^\star(M')$ for each $\Sigma'$-structure $M'$, because, as $U \circ H = U' = U \circ (m_H)^\star$, then for each $\Sigma'$-homomorphism $(M' \overset{g}\to N')$ we will have
$$H(M' \overset{g}\to N') = (m_H)^\star(M' \overset{g}\to N').$$

{\em Claim} $H$ and $(m_H)^\star$ coincide on free $\Sigma'$-structures:\\
Indeed, consider a set $Y$ and the diagram below:

\[\xymatrix{
Y\ar[r]^{\sigma_{Y}}\ar[dr]_{\sigma'_{Y}}&UFY\ar[d]^{U\eta_{H}(Y) =U\eta_{m^{\star}_{H}}(Y)}\\
&U'F'Y
}\]

As $U \circ H = U' = U \circ (m_H)^\star$ and due to the universal property of $\sigma_Y$, them
$$(FY\overset{\eta_{H}(Y)}\to HF'Y) \ = \ (FY\overset{\eta_{m^{\star}_{H}}(Y)}\to m_{H}^{\star}F'Y)$$
 as morphisms of $\Sigma-Str$, hence
$$(+) \ H(F'Y) = m_{H}^{\star}(F'Y).$$

Now we will prove the general case: $H(M') = (m_H)^\star(M')$, for each $M'\in \Sigma'-Str$. Note that $UH(M') = U'(M') = U(m_H)^\star(M')$, thus the $\Sigma$-structures $H(M')$ and $(m_H)^\star(M')$ shares the same underlying set. We must show the the interpretation of all $\Sigma$-symbols in $H(M')$ and $(m_H)^\star(M')$ coincide.

Let $\varepsilon': F'U'\Rightarrow Id_{\Sigma'-Str}$ be the natural transformation that is the  co-unit of the adjunction between $F'$ and $U'$. It is clear that, for each $M' \in \Sigma-Str$, $\varepsilon'_{M'} : F'U'(M') \thra M'$ is a surjective $\Sigma'$-homomorphism, thus the  Isomorphism Theorem gives the following commutative diagram:

\[\xymatrix{
F'U'M'\ar@{->>}[r]^{\varepsilon'_{M'}}\ar@{->>}[dr]_{q_{M'}}& M'\\
&\frac{F'U'M'}{ker(\varepsilon_{M'})}\ar[u]_{\cong \ \bar{q}_{M'}}
}\]

In particular,  the $\Sigma'$-structure $M'$, on the underlying set $U'(M')$, is completely determined by the surjective $\Sigma'$-homomorphims $\varepsilon'_{M'} : F'U'(M') \thra M'$.

Applying $H$ and $m^{\star}_{H}$ to $\varepsilon'_{M'} : F'U'(M') \thra M'$ we obtain the {\em surjective} $\Sigma$-homomorphisms

\[\xymatrix{
HF'U'(M')\ar@{->>}[r]^{H(\varepsilon'_{M'})}\ar@{->>}[dr]_{H(q_{M'})}&H(M')&&m^{\star}_{H}F'U'(M')\ar@{->>}[r]^{m^{\star}_{H}(\varepsilon'_{M'})}\ar@{->>}[dr]_{m^{\star}_{H}(q_{M'})}&m^{\star}_{H}(M')\\
&H(\frac{F'U'M'}{ker(\varepsilon_{M'})})\ar[u]_{\cong \  H(\bar{q}_{M'})}&&& m^{\star}_{H}(\frac{F'U'M'}{ker(\varepsilon_{M'})})\ar[u]_{\cong \ m^{\star}_{H}(\bar{q}_{M'})}
}\]

By $(+)$ above, we have $H(F'(U'(M'))) = m^{\star}_{H}(F'(U'(M')))$, as $\Sigma$-structures. Now, as $U \circ H = U' = U \circ m^{\star}_{H}$, we have

$$ (UHF'U'(M') \overset{UH(\varepsilon'_{M'})}\thra UH(M')) \ = \
 (U'F'U'(M') \overset{U'(\varepsilon'_{M'})}\thra U'(M')) \  = $$
$$ = (Um_H^\star F'U'(M') \overset{Um_H^\star(\varepsilon'_{M'})}\thra Um_H^\star(M')).$$

Thus the $\Sigma$-structures $H(M')$ and $m_H^\star(M')$ on the same underlying set coincide, since they are determined by the same surjective $\Sigma$-homomorphism.





\qed

We will denote the inverse (contravariant) functors in the Theorem above by:

\[\xymatrix{
E_S:\Sf\ar[r]&\Sf^{\dag}&&E_S^{\dag}:\Sf^{\dag}\ar[r]&\Sf\\
\Sigma\ar[d]_{h}&\Sigma-Str&&\Sigma-Str&\Sigma\ar[d]^{m_{H}}\\
\Sigma'&\Sigma'-Str\ar[u]_{h^{\star}}&&\Sigma'-Str\ar[u]_{H}&\Sigma'
}\]

The characterization Theorem \ref{Sf-iso} provides some interesting

\begin{Cor} \label{Sf-cor} Let $H : \Sigma'-Str \to \Sigma-Str$ be a signature functor. Then:

(a) $H$ preserves, {\em strictly}, the following constructions: substructures, products, directed inductive limits,  reduced products, congruences and quotients.

(b) $H$ has a left adjunct $G : \Sigma-Str \to \Sigma'-Str$ with unity of the adjunction  $\lambda : id_{\Sigma-str} \Rightarrow H\circ G$. Moreover $G$ and $\lambda$ can be chosen such that $G \circ F = F'$ and $\lambda_{F(Y)}  = \eta_H(Y) : F(Y) \to H(F'(Y))$, for each set $Y$ and, in particular, from Proposition \ref{h*pres}, for each $Y \sub X$,  $\eta_{H}(Y) = (\check{m}_H)_{\rest Y} : F(Y) \to F'(Y)^{m_H}$.
\end{Cor}

\Dem
(a) This follows from  \ref{h*} and characterization Theorem above.

(b) By characterization Theorem above and Proposition \ref{adjunct-prop}.(a), the functor $H$ has a left adjunct $G$ and, by
Proposition \ref{diag-adj}.(a) $G \circ F \cong F'$. Now we will  analyze the additional  conditions. As adjunct functors are determined up to natural isomorphism by the choice of universal arrows, it is enough to show that, for each set $Y$, the $\Sigma$-homomorphism $\eta_H(Y) : F(Y) \to H(F'(Y))$ is such that for each $M' \in Obj(\Sigma'-Str)$ and each $\Sigma$-homomorphism $f : F(Y) \to H(M')$, there is an unique $\Sigma'$-homomorphism $f' : F'(Y) \to M'$ such that $H(f') \circ \eta_H(Y) = f$. I.e., we must show that, for each $M' \in Obj(\Sigma'-Str)$, the mapping $f' \in \Sigma'-Str(F'(Y),M') \overset{t}\mapsto H(f') \circ \eta_H(Y) \in \Sigma-Str(F(Y), H(M'))$ is a bijection. Consider the bijections given by the pairs of adjunct functors $(F,U)$ and $(F',U')$:

$$f \in \Sigma-Str(F(Y), H(M')) \overset{j}\mapsto U(f) \circ \sigma_Y \in Set(Y, U(H(M')))$$

$$f' \in \Sigma'-Str(F'(Y),M') \overset{j'}\mapsto U'(f') \circ \sigma'_Y \in Set(Y, U'(M'))$$

As $Set(Y, U'(M'))  = Set(Y, U(H(M')))$ and $U(\eta_H(Y)) \circ \sigma_Y = \sigma'_Y$ we conclude that $j \circ t = j'$, i.e., the diagram below commutes

\[\xymatrix{
Set(Y, U'(M'))\ar[r]^{=}& Set(Y, U(H(M'))\\
\Sigma'-Str(F'(Y),M')\ar[u]^{\cong \ j'}\ar[r]_{t}&\Sigma-str(F(Y), H(M'))\ar[u]_{j \ \cong}
}\]

Thus, as $j$ and $j'$ are bijections, then $t$ is a bijection. This entails the additional results.
\qed

\vcinco

Now, having a detailed functorial encoding of (flexible) signature morphisms, we can proceed to a functorial description of logical morphisms between algebraizable logics.

\vcinco

\begin{Lem} \label{bareta-le}
Let $I : \cK \hookr \Sigma-Str$ and $I' : \cK' \hookr \Sigma'-Str$  full inclusions, where $\cK$ and $\cK'$ are quasivarieties. Let $H : \Sigma'-Str \ra \Sigma-Str$ be a  signature functor such that it restricts (uniquely) to a functor $H\!\!\rest : \cK' \to \cK$ (thus $I \circ H\!\!\rest = H \circ I'$).  Keeping the notation in Remark \ref{adjQV-re}, for each set $Y$, let (by the universal property of $t_Y$) $\bar{\eta}_H(Y) : LF(Y) \ra H\!\!\rest(L'F'(Y))$ be the unique $\cK$-morphism such that\ $(Y \overset{t_Y}\to UILF(Y)  \overset{UI(\bar{\eta}_Y)}\to UIH\!\!\rest L' F'(Y)) \ =$ \\ $(Y \overset{t'_Y}\to U'I'L'F'(Y))$. Then:

(a) $(\bar{\eta}_H(Y))_{Y \in Set}$ is a natural transformation $\bar{\eta}_{H}: L \circ F \ra H\!\!\rest \circ L'\circ F'$.

(b) Both the diagrams below commute

\[\xymatrix{
Y\ar[r]^{\sigma_{Y}}\ar@/^2pc/[rr]^{t_{Y}}\ar[d]_{id_Y}&UF(Y)\ar[r]^{U(q_{F(Y)})}\ar[d]_{U(\eta_{H}(Y))}&UILF(Y)\ar[d]^{UI(\bar{\eta}_H(Y))}\\
Y\ar[r]_{\sigma'_{Y}}\ar@/_2pc/[rr]_{t'_{Y}}&UHF'(Y)\ar[r]_{UH(q'_{F'(Y)})}&UIH\!\!\rest L'F'(Y)
}\]
\vcinco

\[\xymatrix{
 F(Y)\ar[dd]_{\eta_H(Y)}\ar[rr]^{q_{F(Y)}}&
&ILF(Y)\ar[dd]^{I(\bar{\eta}_H(Y))}\\
&&&\\
 HF'(Y)\ar[rr]_{H(q'_{F'(Y)})}& &IH\!\!\rest L'F'(Y)
}\]

(c) $H$ and $H\!\!\rest$ have  left adjuncts, respectively  $G : \Sigma-Str \to \Sigma'-Str$  and  $\bar{G} : \cK \to \cK'$, the respective unities of the adjunctions  $\lambda : id_{\Sigma-str} \Rightarrow H\circ G$ and $\bar{\lambda} : id_{\cK} \Rightarrow H\!\rest\circ \bar{G}$. Moreover $G, \bar{G}$ and $\lambda, \bar{\lambda}$ can be chosen such that:\\
$\bullet$ \ $G \circ F = F'$ and $\bar{G} \circ L \circ F = L' \circ F' = L' \circ G \circ F$;\\
$\bullet$ \  $\lambda_{F(Y)}  = \eta_H(Y) : F(Y) \to H(F'(Y))$ and $\bar{\lambda}_{LF(Y)}  = \bar{\eta}_H(Y) : LF(Y) \to$ \\
$H\!\!\rest(L'F'(Y))$, for each set $Y$.

\end{Lem}

\Dem Item (a) follows in an analogous fashion to the proof of Lemma \ref{sigfunt-le}.(a): by analyzing the commutativity of the diagram below from the universal property of $t_Y$,  for each function $f : Y \to Z$.

\[\xymatrix{
Y\ar[r]^{t_{Y}}\ar@/^2pc/[rr]^{t'_{Y}}\ar[d]_{f}&UILF(Y)\ar[r]^{UI(\bar{\eta}_{H}(Y))}\ar[d]_{UILF(f)}&UIH\!\rest L' F'(Y)\ar[d]^{UIH\!\rest L' F'(f)}\\
Z\ar[r]_{t_{Z}}\ar@/_2pc/[rr]_{t'_{Z}}&UILF(Z)\ar[r]_{UI(\bar{\eta}_{H}(Z))}&UIH\!\rest L'F'(Z)
}\]
\vcinco

Item (b) follows in an analogous fashion to the proof of Lemma \ref{sigfunt-le}.(a): the top diagram commutes, by analyzing the commutativity of the diagram below from the universal property of $\sigma_Y$; the bottom diagram commutes since the functor $U$ is faithful and the inner right square in  the top diagram commutes.

Item (c) follows in an analogous fashion to the proof of Corollary \ref{Sf-cor}.(b): first, by applying Proposition \ref{diag-adj},  and then, by a diagram chase to shows that, for each $M' \in \cK'$, the mapping $f' \in \cK'(L'F'(Y),M') \mapsto H\!\rest(f') \circ \bar{\eta}_H(Y) \in \cK(LF(Y), H\!\rest(M'))$ is a bijection.

\qed

\begin{Prop}\label{QV-dense-prop-2} Let  $l = (\Sigma, \vdash) \in \Lf$ and $a, a' \in \Af$.

(a) Let  $h : l \ra a'$ be a $\Lf$-morphism. Then are equivalent:

(a1) \  $h$ is a $\Delta$-dense $\Lf$-morphism.

(a2) \ The functor $h^\star\!\!\rest : QV(a') \ra \Sigma-Str$ is full, faithful, injective on objects and satisfies the heredity  condition (see \ref{adjunct-prop}.(b4)).

(b) Let  $h: a \ra a'$ be a $\Af$-morphism. Then are equivalent:

(b1) \  $h$ is a $\Delta$-dense $\Af$-morphism.

(b2) \ The functor $h^\star\!\!\rest : QV(a') \ra QV(a)$ is full, faithful, injective on objects and satisfies the heredity  condition.

\end{Prop}

\Dem The implications (a1) $\Rightarrow$ (a2) and (b1) $\Rightarrow$ (b2) were established in Proposition \ref{QV-dense-prop}.

\underline{(a1) $\Rightarrow$ (a2)}:
by Theorem \ref{Sf-iso},  Lemma \ref{bareta-le}.(b), Remark \ref{freerest} and Corollary \ref{Sf-cor}.(b), the following diagram commutes, for each $Y \sub X$.

\[\xymatrix{
 F(Y)\ar[dd]_{\check{h}\!\!_{\rest Y}}\ar@{->>}[rr]^{id_{F(Y)}}&
&F(Y)\ar[dd]^{I(\bar{\eta}_{h^\star}(Y))}\\
&&&\\
 h^\star(F'(Y))\ar@{->>}[rr]_{h^\star(q'_{F'(Y)})}& &I'(h^\star\!\rest(F'(Y)/\Delta'\!\rest))
}\]

By  hypothesis (a1),  Lemma \ref{bareta-le}.(c) and Proposition \ref{adjunct-prop}.(b), the $\Sigma'$- homomorphism $\bar{\eta}_H(Y) : F(Y) \to I'(F'(Y)/\Delta'\!\rest)$ is {\em surjective}. Thus a diagram chase shows that for each $\phi' \in F'(Y')$ there is $\phi \in F(Y)$ such that $\vdash' \check{h}(\phi) \Delta' \phi'$. Therefore, the $\Lf$-morphism $h : l \to a$ is $\Delta$-dense.

\underline{(b1) $\Rightarrow$ (b2)}: is proved in an analogous way, by a chase on the commutative diagram below

\[\xymatrix{
 F(Y)\ar[dd]_{\check{h}\!\!_{\rest Y}}\ar@{->>}[rr]^{q_{F(Y)}}&
&I(F(Y)/\Delta\!\rest)\ar@{->>}[dd]^{I(\bar{\eta}_{h^\star}(Y))}\\
&&&\\
 h^\star(F'(Y))\ar@{->>}[rr]_{h^\star(q'_{F'(Y)})}& &I'(h^\star\!\rest(F'(Y)/\Delta'\!\rest))
}\]

\qed

As density and $\Delta$-density of morphisms coincide on Lindenbaum algebraizable logics, we immediately obtain the

\begin{Cor}\label{QV-dense-cor-2} Let  $l = (\Sigma, \vdash) \in \Lf$ and $a, a' \in \cA^c_f$.

(a) Let  $h : l \ra a'$ be a $\Lf$-morphism. Then are equivalent:

(a1) \  $h$ is a dense $\Lf$-morphism.

(a2) \ The functor $h^\star\!\!\rest : QV(a') \ra \Sigma-Str$ is full, faithful, injective on objects and satisfies the heredity  condition .

(b) Let  $h: a \ra a'$ be a $\cA_f^c$-morphism. Then are equivalent:

(b1) \  $h$ is a dense $\cA_f^c$-morphism.

(b2) \ The functor $h^\star\!\!\rest : QV(a') \ra QV(a)$ is full, faithful, injective on objects and satisfies the heredity  condition.

\end{Cor}

Having in mind the Definitions \ref{sigfunct} and \ref{Sf+}, it is natural to consider the following

\begin{Df} \label{Lf+}

(a) Let $a = (\Sigma, \vdash), a'= (\Sigma', \vdash')$ be algebraizable logics.  A functor $H : \Sigma'-Str \ra \Sigma-Str$ will be called a  "BP-functor", $H$ is a signature functor also satisfying $(l1), (l2), (l3)$: \\
$(l1)$ $H$ has a (unique) restriction to the associated  quasivarieties $H\!\!\rest : QV(a') \ra QV(a)$;\\
There are algebraizing pairs $(\Delta,(\delta,\varepsilon))$ and $(\Delta',(\delta',\varepsilon'))$ of, respectively, $a$ and $a'$ such that:\\
$(l2)$ $\check{m}_{H}(\Delta)\dashv ' \vdash\Delta'$;\\
$(l3)$ $\check{m}_{H}(\delta) \equiv \check{m}_H(\varepsilon) =\!|_{QV(a')}|\!= \delta' \equiv \varepsilon'$.

It is straightforward that:\\
$\bullet$ $id_{\Sigma-Str} : \Sigma-Str \ra \Sigma-Str$ is a BP-functor;\\
$\bullet$ If $(\Sigma-Str \overset{H}\leftarrow \Sigma'-Str \overset{H'}\leftarrow
\Sigma''-Str)$ are BP-functors, then $H\circ H' : \Sigma''-Str \to \Sigma-str$ is a BP-functor.\footnote{Note that, for each $M'' \in QV(a'')$, $(M'')^{m_{H'}} = H'(M'') \in QV(a')$ and $(M'')^{m_{H'}} \vDash_{\Sigma'} \check{m}_{H}(\delta) \equiv \check{m}_H(\varepsilon) \leftrightarrow \delta' \equiv \varepsilon'$ iff $M'' \vDash_{\Sigma''} \check{m}_{H'}(\check{m}_{H}(\delta)) \equiv \check{m}_{H'}(\check{m}_H(\varepsilon)) \leftrightarrow \check{m}_{H'}(\delta') \equiv  \check{m}_{H'}(\varepsilon')$.}
\\

(b) Denote $\cA_f^\dagger$ the category with:\\
$\bullet$  {\em Objects:} are pairs $(\Sigma-Str, a)$ where $a = (\Sigma, \vdash)$ is an algebraizable logic;\\
$\bullet$ {\em Arrows:} are BP-functors $(\Sigma'-Str, a') \overset{H}\to (\Sigma-Str, a)$;\\
$\bullet$ {\em identities and composition:} as (BP-)functors.
\\

(c) Denote $Lind(\cA_f)^\dagger$ the full subcategory of $\cA_f^\dagger$ with objects, the  pairs $(\Sigma-Str, a)$ where $a = (\Sigma, \vdash)$ is a Lindenbaum algebraizable logic.
\end{Df}
\vtres

Below we present the results that encompass most part of the present work

\begin{Teo} \label{Af-iso}
The pair of inverse anti-isomorphisms of categories $\Sf \overset{E_S}{\underset{E_S^\dagger}{\rightleftarrows}} \Sf^\dagger$ in Theorem \ref{Sf-iso}  "restricts", via the forgetful functors $\cA_f \ra \Sf$ and $\cA_f^\dagger \ra \Sf^\dagger$,  to a pair of inverse anti-isomorphisms of categories $\cA_f \overset{E_A}{\underset{E_A^\dagger}{\rightleftarrows}} \cA_f^\dagger$.

\[\xymatrix{
E_A:\Af\ar[r]&\Af^{\dag}&&E_A^{\dag}:\Af^{\dag}\ar[r]&\Af\\
a=(\Sigma, \vdash)\ar[d]_{h}&(\Sigma-Str, a)&&(\Sigma-Str,a)&a\ar[d]^{m_{H}}\\
a'=(\Sigma',\vdash')&(\Sigma'-Str,a')\ar[u]_{h^{\star}}&&(\Sigma'-Str,a')\ar[u]_{H}&a'
}\]

\[\xymatrix{
\cA_f\ar[dd]_{Forget}\ar[rr]_{E_A}&
&\cA_f^\dagger\ar[dd]^{Forget}\ar[ll]_{E_A^\dagger}\\
&&&\\
 \Sf\ar[rr]_{E_S}& &\Sf^\dagger\ar[ll]_{E_S^\dagger}
}\]

Moreover, if $h \in \cA_f(a,a')$ and $H \in \cA_f^\dagger((\Sigma'-Str, a'), (\Sigma-Str, a))$ are in correspondence, then  the  pair of inverse anti-isomorphisms $(E_A, E_A^\dagger)$ is such that:

(a) It establishes a correspondence between the equivalence class $\{h' \in \cA_f(a,a'): [h]_{\approx}= [h']_{\approx}  \in \overline{\cA_f}(a,a')\}$ and   the   equivalence class $\{H' \in \cA_f^\dagger((\Sigma'-Str, a'), (\Sigma-Str, a)) : H'\rest = H\rest\}$.

(b)  $[h]_{\approx}$ is a $\overline{\cA_f}$-isomorphism  \hem $\Leftrightarrow$ \hem $H\!\!\rest$ is an  isomorphism between quasivarieties.

(c)  $h$ is a $\Delta$-dense morphism \hem $\Leftrightarrow$ \hem  $H\!\!\rest$ is full, faitful, injective on object and heredity.

\end{Teo}

\Dem After the pair of ("restricted") inverse anti-isomorphisms $(E_A, E^\dagger_A)$ were established, then: item (a)  follows from Proposition \ref{QV-quo-prop}.(b); item (b) follows from Proposition \ref{QV-iso-prop}; item (c) follows from Proposition  \ref{QV-dense-prop-2}.(b).\\

It follows from  directly from Theorem \ref{Sf-iso} and the definitions of the object part of the functors $(E_A, E_A^\dagger)$ that they establishes an well defined pair of inverse bijections between the classes of objects $Obj(\cA_f)$ and $Obj(\cA_f^\dagger)$.\\

If we establish that the (arrow) mappings below are well defined:\\
$h \in \cA_f(a,a') \overset{E_A}\mapsto h^\star \in \cA_f^\dagger((\Sigma'-Str,a'),(\Sigma-Str, a))$;\\
$H \in \cA_f^\dagger((\Sigma'-Str,a'),(\Sigma-Str, a)) \overset{E_A^\dagger}\mapsto m_H \in \cA_f(a,a'),$\\
then it will follow from Theorem \ref{Sf-iso} that
the pair of inverse anti-isomorphisms $\Sf \overset{E_S}{\underset{E_S^\dagger}{\rightleftarrows}} \Sf^\dagger$ in Theorem \ref{Sf-iso}  "restricts"  to a pair of inverse anti-isomorphisms $\cA_f \overset{E_A}{\underset{E_A^\dagger}{\rightleftarrows}} \cA_f^\dagger$.\\

Let $h \in \cA_f(a,a')$. By Proposition \ref{h*pres},   $h^\star : \Sigma'-str \to \Sigma-Str$ is a signature functor  and, by  Proposition \ref{restrict}, it restricts (uniquely) to a functor $h^{\star}\!\!\rest:QV(a')\to QV(a)$: thus condition $(l1)$ is fulfilled. By Theorem \ref{Sf-iso}, $m_{h^{\star}}= h$; as $h$ preserves algebraizable pairs, then Fact \ref{uniq}.(a)  ensures that the conditions $(l2)$ and $(l3)$ are satisfied. Therefore $E_A$ is an well defined functor.\\

Let $H \in \cA_f^\dagger((\Sigma'-Str,a'),(\Sigma-Str, a))$. Lemma \ref{sigfunt-le}.(c) entails that $m_H : \Sigma \to \Sigma'$ is a $\Sf$-morphism. Conditions $(l2)$ and $(l3)$ and Fact \ref{uniq}.(b) ensures that $m_H$ preserves algebraizing pairs.  It remains to show that $m_H$ is a $\Lf$-morphism, i.e. given $\Gamma\cup\{\varphi\}\subseteq F(X)$, we must have
\[\Gamma\vdash \varphi\ \Rightarrow\ \check{m}_H[\Gamma]\vdash' \check{m}_H(\varphi)\]
But, as $a$ and $a'$ are algebraizable logics, it is enough to prove that
$$\{\varepsilon(\psi)\equiv\delta(\psi);\ \psi\in\Gamma\}\models_{QV(a)}\varepsilon(\varphi)\equiv\delta(\varphi)\ \Rightarrow$$
$$\{\varepsilon'(\check{m_{H}}(\psi))\equiv\delta'(\check{m_{H}}(\psi));\ \psi\in\Gamma\}\models_{QV(a')}\varepsilon'(\check{m_{H}}(\varphi))\equiv\delta'(\check{m_{H}}(\varphi)).$$

Let $M'\in QV(a')$ and suppose that $M' \models_{\Sigma'} \varepsilon'(\check{m_{H}}(\psi))\equiv\delta'(\check{m_{H}}(\psi))$ for each $\psi\in\Gamma$. As $m_{H}$ satisfies condition $(l3)$, then holds, for each $\psi\in\Gamma$,
\[M' \models_{\Sigma'} \check{m_{H}}(\varepsilon)(\check{m_{H}}(\psi))\equiv\check{m_{H}}(\delta)(\check{m_{H}}(\psi))\]
I.e.:
\[M' \models_{\Sigma'} \check{m_{H}}(\varepsilon(\psi))\equiv\check{m_{H}}(\delta(\psi))\]

 By Theorem \ref{Sf-iso}, $H = (m_H)^\star$, thus we get

\[H(M') \models_{\Sigma} \varepsilon(\psi)\equiv\delta(\psi)\]

From the hypothesis, $H(M') \in QV(a)$, and as $\{\varepsilon(\psi)\equiv\delta(\psi);\ \psi\in\Gamma\}\models_{QV(a)}\varepsilon(\varphi)\equiv\delta(\varphi)$, we obtain

\[H(M') \models_{\Sigma} \varepsilon(\varphi)\equiv\delta(\varphi)\]

Therefore, as above,

\[M' \models_{\Sigma'} \check{m_{H}}(\varepsilon(\varphi))\equiv\check{m_{H}}(\delta(\varphi))\]

and

 $$M' \models_{\Sigma'} \varepsilon'(\check{m_{H}}(\varphi))\equiv\delta'(\check{m_{H}}(\varphi)).$$

As $M' \in QV(a')$ was taken arbitrarily, then $\{\varepsilon'(\check{m_{H}}(\psi))\equiv\delta'(\check{m_{H}}(\psi));\ \psi\in\Gamma\}\models_{ QV(a')}\varepsilon'(\check{m_{H}}(\varphi))\equiv\delta'(\check{m_{H}}(\varphi))$.

Summing up, $m_{H}$ is a logical morphism that preserves algebraizable pairs. Therefore $E_A^\dagger$ is an well defined functor. This finishes the proof.
\qed
\vtres

Restricting the result above to the setting of Lindenbaum algebraizable logics, we obtain the

\begin{Cor} \label{LindAf-iso}
The pair of inverse anti-isomorphisms of categories  $\Af \overset{E_A}{\underset{E_A^\dagger}{\rightleftarrows}} \Af^\dagger$ in Theorem \ref{Af-iso}  "restricts", via the (full) inclusion functors $Lind(\cA_f) \hookr \Af$ and $Lind(\cA_f)^\dagger \hookr \Af^\dagger$,  to a pair of inverse anti-isomorphisms of categories  $Lind(\cA_f) \overset{E_L}{\underset{E_L^\dagger}{\rightleftarrows}} Lind(\cA_f)^\dagger$.

\[\xymatrix{
Lind(\cA_f)\ar[dd]_{Incl}\ar[rr]_{E_L}&
&Lind(\cA_f)^\dagger\ar[dd]^{Incl}\ar[ll]_{E_L^\dagger}\\
&&&\\
 \Af\ar[rr]_{E_A}& &\Af^\dagger\ar[ll]_{E_A^\dagger}
}\]

Moreover, if $h \in Lind(\cA_f)(a,a')$ and $H \in Lind(\cA_f)^\dagger((\Sigma'-Str, a'), (\Sigma-Str, a))$ are in correspondence, then  the  pair of inverse anti-isomorphisms $(E_L, E_L^\dagger)$ is such that:

(a) It establishes a correspondence between the equivalence class \[\{h' \in Lind(\cA_f)(a,a'): [h]_{\dashv \ \vdash}= [h']_{\dashv \ \vdash}  \in QLind(\cA_f)(a,a')\}\]
and   the   equivalence class \[\{H' \in Lind(\cA_f)^\dagger((\Sigma'-Str, a'), (\Sigma-Str, a)) : H'\rest = H\rest\}.\]

(b)  $[h]_{\dashv \ \vdash}$ is a $QLind(\cA_f)$-isomorphism  \hem $\Leftrightarrow$ \hem $H\!\!\rest$ is an  isomorphism between quasivarieties.

(c)  $h$ is a dense morphism \hem $\Leftrightarrow$ \hem  $H\!\!\rest$ is full, faitful, injective on object and heredity.

\end{Cor}

\Dem It is clear that $(E_A, E^\dagger_A)$ establishes a bijective correspondence between the subclasses $Obj(Lind(\cA_f))$ and $Obj(Lind(\cA_f)^\dagger)$. As $Lind(\cA_f) \hookr \Af$ and $Lind(\cA_f)^\dagger \hookr \Af^\dagger$ are full subcategories, then $(E_A, E^\dagger_A)$ restricts to a pair of inverse anti-isomorphisms $Lind(\cA_f) \overset{E_L}{\underset{E_L^\dagger}{\rightleftarrows}} Lind(\cA_f)^\dagger$.

On the additional results: item (a)  follows from Corollary \ref{QV-quo-cor}.(b); item (b) follows from Corollary \ref{QV-iso-cor}; item (c) follows from Corollary  \ref{QV-dense-cor-2}.(b).
\vcinco
\qed

\section{Future works}

We believe that the notions and results here presented  can play the role of initial steps towards  a development of a representation theory of (propositional) logics through the category theory. So we intend apply  mathematical devices to obtain information about logics and the study meta-logical properties. In particular, we want apply this representation theory to more involved logics as the  Logics of Formal Inconsistency --LFIs--(\cite{BCC1}, \cite{BCC2}).

Influenced  by the processes  of analysis and synthesis of logics in the combining of logics field, we have considered in  \cite{MaPi1} some notions of left/right Morita equivalence of logics (notions weaker than isomorphism of logics), based on the functorial encoding of logical morphisms.
We intend also obtain information on a logic through  comparison with other "well behaved" of Lindenbaum-algebraizable logics: a kind of "local-global principle" approach.

In \cite{MaPi2}, we consider another application of the ideas here developed  to investigate a functorial treatment of "G\"odel translations" in the setting of (two) Lindenbaum algebraizable logics, as we already mentioned in \ref{godel}: there we define an institution to each "kind" of Lindenbaum algebraizable logic and derive a Glivenko's theorem  between two "kinds" of such logics, by defining a morphism between its associated institutions.

Inspired by results on algebraizable logics  that identifies the occurrence of a certain meta-logical property  by a property on the quasivariety associated to the algebraizable logic\footnote{For instance, Beth property \cite{Ho} and Craig interpolation property \cite{Cze2} for consequence relation (see \cite{FJP} for a compilation of results).}, we begin in \cite{AJMP} an attempt of a categorial local-global  analysis of meta-logical properties of a given logic through its diagram of (Lindenbaum) algebraizable logics.
We intend analyze the behavior  under categorial constructions as products and directed colimits (among others), of  this  local-global approach to meta-logic properties of  logics .


\end{document}